\documentclass[english,12pt,a4paper,twoside,reqno]{article}
\usepackage{amsfonts}
\usepackage{amsmath}
\usepackage{stmaryrd}
\usepackage{amssymb}
\usepackage{mathrsfs}
\usepackage[T1]{fontenc}
\usepackage{graphics}
\usepackage{color}
\usepackage[english]{babel}
\usepackage{bm}
\usepackage{enumerate}
\usepackage{eufrak}
\usepackage{euscript}
\numberwithin{equation}{section}
\usepackage{newlfont}
\usepackage{latexsym}
\usepackage{graphicx}
\usepackage{a4wide}

\usepackage{cite}        % Sorted citations, if these are numerical
\usepackage{url,float}         % Define \url{http://...} command

\title{Evolution equations in discrete and continuous time for nonexpansive operators in  Banach spaces}
\author{Guillaume Vigeral\\
\footnotesize{Equipe Combinatoire et Optimisation, CNRS FRE3232, Université P. et M. Curie, Paris 6, UFR 929}\\
 \footnotesize{175 rue du Chevaleret, 75013 Paris, France}\\
 \footnotesize{guillaumevigeral@gmail.com}}
\date{\today}
\newtheorem{theorem}{Theorem}[section]
\newtheorem{proposition}[theorem]{Proposition}
\newtheorem{lem}[theorem]{Lemma}
\newtheorem{corollary}[theorem]{Corollary}

\newtheorem{déf}[theorem]{Definition}

\newtheorem{remark}[theorem]{Remark}
\newtheorem{example}[theorem]{Example}
\newenvironment{proof}[1][Proof]{\noindent \textbf{#1.}~ }
{\hfill\rule{2mm}{2mm} \vspace{\parskip} }

\begin{document}
\maketitle
 \begin{abstract}
%We consider some evolution equations in a Banach space involving a
%nonexpansive operator. Our motivation comes from the study of
%zero-sum two players repeated games, whose values satisfy some discrete dynamics. \\
%We prove that the solutions of those continuous-time equations and
%discrete dynamics have the same asymptotic behavior, even in the
%absence of convergence. In the process we also establish a new
%exponential formula for explicit schemes, as well as a Kobayashi
%inequality.
We consider some discrete and continuous dynamics in a Banach space
involving a non expansive operator $J$ and a corresponding family of
strictly contracting operators $\Phi(\lambda,x):=\lambda
J(\frac{1-\lambda}{\lambda}x)$ for $\lambda\in]0,1]$. Our motivation
comes from the study of two-player zero-sum repeated games, where
the value of the $n$-stage game (resp. the value of the
$\lambda$-discounted game) satisfies the relation
$v_n=\Phi(\frac{1}{n},v_{n-1})$ (resp.
$v_\lambda=\Phi(\lambda,v_\lambda)$) where $J$ is the Shapley
operator of the game. We study the evolution equation
$u'(t)=J(u(t))-u(t)$ as well as associated Eulerian schemes,
establishing a new exponential formula and a Kobayashi-like
inequality for such trajectories. We prove that the solution of the
non-autonomous evolution equation
$u'(t)=\Phi(\bm{\lambda}(t),u(t))-u(t)$ has the same asymptotic
behavior (even when it diverges) as the sequence $v_n$ (resp. as the
family $v_\lambda$) when $\bm{\lambda}(t)=1/t$ (resp. when
$\bm{\lambda}(t)$ converges slowly enough to 0).
\end{abstract}

\section{Introduction}
The topic of the asymptotic behavior of trajectories defined through
nonexpansive mappings in Banach spaces arise in numerous domains
such as nonlinear semigroups theory \cite{Ba,Br,CrLi,Ka,Ko,MiOh,Re},
game theory \cite{KoNe,Ne,NeSo,RoSo,So,So2} as well as in discrete
events systems \cite{GaGu,Gu,GuKe}.

Given a nonexpansive function $J$ from a Banach space $X$ to itself,
evolution equation
\begin{equation}
\label{eqintro1}U'(t)=J(U(t))-U(t)
\end{equation}
\noindent is a particular case of
the widely-studied
\begin{equation*}
U'(t)\in-A(U(t))
\end{equation*}
\noindent for a maximal monotone operator $A$. Typically, the study
of the asymptotics for such evolution equation and its Eulerian and
proximal discretizations has been made in Hilbert spaces\cite{Br} or
at least assuming some geometric properties in the case of Banach
spaces\cite{KoNe,Re}. Another usual assumption is the non emptiness
of the set $A^{-1}(0)$.

On the other hand, in the framework of two-person zero-sum games
repeated in discrete time, the values $v_n$ and $v_\lambda$ of the
$n$-stage (resp. $\lambda$-discounted) game satisfy respectively:
\begin{equation}\label{eqintrovn}v_n=\frac{J^n(0)}{n}=\Phi\left(\frac{1}{n},v_{n-1}\right)\end{equation}
\begin{equation}\label{eqintrovlambda}v_\lambda=
%\lambda J\left(\frac{1-\lambda}{\lambda}v_\lambda\right)
\Phi(\lambda,v_\lambda)
\end{equation}
\noindent where $J$ is the so-called Shapley operator of the game
and $\Phi(\lambda,x):=\lambda
J\left(\frac{1-\lambda}{\lambda}x\right)$. This operator $J$ is
nonexpansive for the uniform norm, hence $A=I-J$ is a maximal
monotone operator in the sense of \cite{Ka}. However two unusual
facts appears in the study of the asymptotics of those values: first
$A^{-1}(0)$, the set of fixed points of $J$, is generally empty.
Another difficulty lies in the lack of smoothness of the unit ball
$\displaystyle{\mathcal{B}_{\|\cdot\|_{\infty}}}$, which might
induce oscillations of the discrete trajectories defined
above\cite{KoNe}.
%However in that case the set $A^{-1}(0)$ is usually empty, and
%the underlying Banach space lacks any nice geometrical property.

The purpose of this paper is to investigate the relation between
several discrete and continuous dynamics in Banach spaces. Because
our motivation comes from this game-theoretic framework, we neither
make any geometrical assumptions on the unit ball, nor suppose non
emptiness of $A^{-1}(0)$. In continuous time, dynamics that we will
consider are (\ref{eqintro1}) as well as non autonomous evolution
equations of the form
\begin{equation}
\label{eqintro2} u'(t)=\Phi(\bm{\lambda}(t),u(t))-u(t)
\end{equation}
\noindent for some parametrizations $\bm{\lambda}$. We establish
that the quantities defined in (\ref{eqintrovn}) and
(\ref{eqintrovlambda}) behave asymptotically as the solutions of
these various evolution equations. Surprisingly this is true not
only when there is convergence; even when they oscillate we prove
that discrete and continuous trajectories remain asymptotically
close.

Section 2 is devoted to definitions and basic results. In Section 3
we study the relation between the solution $U$ of evolution equation
(\ref{eqintro1}) and related Eulerian schemes, establishing in
particular that $\|v_n-\frac{U(n)}{n}\|$ converges to 0. In the
process we prove that some classical results (e.g.\ exponential
formula\cite{CrLi}, Kobayashi inequality\cite{Ko}) involving the
proximal trajectories for a maximal monotone operator $A$ have an
Eulerian explicit counterpart in the case $A=I-J$. In Section 4 we
consider the non autonomous equation (\ref{eqintro2}). We show that
for $\bm{\lambda}(t)=\frac{1}{t}$ the solution behave asymptotically
as the sequence $v_n$, and that when $\bm{\lambda}$ converges slowly
enough to 0 the solution behave asymptotically as the family
$v_\lambda$.

\section{Discrete time model}

\subsection{Nonexpansive operators}

Let $(X,\|\cdot\|)$ be a Banach space, and $J$ a nonexpansive
mapping from $X$ into itself :
\[ \|J(x)-J(y)\|\leq\|x-y\| \quad \forall(x,y)\in X^2.\]

\noindent We define, for $n\in\mathbb{N}$ and $\lambda\in]0,1]$,
\begin{eqnarray}
\label{defVn}V_n&=&J(V_{n-1})=J^n(0)\\
\label{defVlambda}V_\lambda&=&J((1-\lambda)V_\lambda)
\end{eqnarray}
Notice that $V_\lambda$ is well-defined because
$J((1-\lambda)\cdot)$ is strictly contracting, hence has a unique
fixed point.
\begin{example}\label{exempletrivial}
For any $c\in\mathbb{R}$, the mapping $J$ from $\mathbb{R}$ to
itself defined by $J(x)=x+c$ is nonexpansive. In that case, $V_n=nc$
and $V_\lambda=\frac{c}{\lambda}$.
\end{example}
\noindent These quantities being unbounded in general (see above),
we also introduce their normalized versions
\begin{eqnarray}
v_n&=&\frac{V_n}{n}\\
v_\lambda&=&\lambda V_\lambda
\end{eqnarray}
In the previous example, one gets $v_n=v_\lambda=c$ for all $n$ and
$\lambda$. In general it is easy to prove that these normalized
quantities are bounded:
\begin{lem}\label{vnlambdabornes}
For any $n\in\mathbb{N}$ and $\lambda\in]0,1]$,
\begin{eqnarray}
\|v_n\|&\leq&\|J(0)\| \\
\|v_\lambda\|&\leq&\|J(0)\|.
\end{eqnarray}
\end{lem}

\begin{proof}
Since $J$ is non expansive,
\[\|V_n-V_{n-1}\|=\|J(V_{n-1})-J(V_{n-2})\|\leq\|V_{n-1}-V_{n-2}\|.\]
By induction this implies that \[\|V_n\|\leq n\|V_1\|=n\|J(0)\|.\]
On the other hand, again using the fact that $J$ is non expansive,
\begin{eqnarray*}
\|V_\lambda\|-\|J(0)\|&\leq&\|V_\lambda-J(0)\| \\
&=&\|J((1-\lambda)V_\lambda)-J(0)\|\\
&\leq& (1-\lambda)\|V_\lambda\|\end{eqnarray*} and so
\[\|v_\lambda\|=\lambda\|V_\lambda\|\leq\|J(0)\|.\]
\end{proof}

To underline the link between the families $
\{v_n\}_{n\in\mathbb{N}}$ and $\{v_\lambda\}_{\lambda\in ]0,1]}$ it
is also of interest to introduce the family of strictly contracting
operators $\Phi(\lambda,\cdot)$, $\lambda\in]0,1]$, defined by

\begin{equation}\label{defPhi}
\Phi(\lambda,x)=\lambda J\left(\frac{1-\lambda}{\lambda}x\right).
\end{equation}
The function $\Phi(\lambda,\cdot)$ can be seen as a perturbed
recession function of $J$: because of the nonexpansiveness of $J$,
\begin{equation}
\lim_{\lambda\rightarrow 0}\Phi(\lambda,x)=\lim_{\lambda\rightarrow
0}\lambda J\left(\frac{x}{\lambda}\right)=\lim_{t\rightarrow
+\infty} \frac{J\left(tx\right)}{t}
\end{equation}
which is the definition of the recession function of $J$\cite{Roc}.

The quantities $v_n$ and $v_\lambda$ then satisfy the relations
\begin{eqnarray}
\label{defvn}v_{n}&=&\Phi\left(\frac{1}{n},v_{n-1}\right)\ ;\ v_0=0\\
\label{defvlambda}v_\lambda&=&\Phi(\lambda,v_\lambda)
\end{eqnarray}
Notice that since $\Phi(\lambda,\cdot)$ is strictly contracting, any
sequence $w_n\in X$ satisfying

\begin{equation}\label{eqrecvlambda}
 w_n=\Phi(\lambda,w_{n-1})
\end{equation}
converges strongly to $v_\lambda$ as $n$ goes to $+\infty$.

\subsection{Shapley operators}
 An important application, which is our main motivation, is
obtained in the framework of zero-sum two player repeated
games\cite{So2}. For example take the simple case of a stochastic
game with a finite state space $\Omega$, compact move sets $U$ and
$V$ for player 1 and 2 respectively, payoff $g$ from $U\times
V\times \Omega$ to $\mathbb{R}$, and transition probability $\rho$
from $U\times V\times \Omega$ to $\Delta(\Omega)$ (the set of
probabilities on $\Omega$).  Let $S=\Delta_f(U)$ (resp.
$T=\Delta_f(V)$) the sets of probabilities on U (resp. V) with
finite support; we still denote by $g$ and $\rho$ the multilinear
extensions from $U\times V$ to $S\times T$ of the corresponding
functions.

The game is played as follow: an initial stage $\omega_1\in\Omega$
is given, known by each player. At each stage $m$, knowing past
history and current state $\omega_m$, player 1 (resp. player 2)
chooses $\sigma\in S$ (resp. $\tau\in T$). A move $a_m$ of player 1
(resp. $b_m$ of player 2) is drawn accordingly to $\sigma$ (resp.
$\tau$). The payoff $g_m$ at stage $m$ is then $g(a_m,b_m,\omega_m)$
and $\omega_{m+1}$, the state at stage $m+1$, is drawn accordingly
to $\rho(a_m,b_m,\omega_m)$.

There are several ways of evaluating a payoff for a given infinite
history:
\begin{itemize}
\item $\frac{1}{n}\sum_{m=1}^n g_m$ is the payoff of the $n-$stage game
\item $\lambda \sum_{m=1}^{+\infty} (1-\lambda)^{i-1} g_m$ is the payoff of the $\lambda-$discounted game.
\end{itemize}
For a given initial state $\omega$, we
denote the values of those games by $v_n(\omega)$ and $v_{\lambda}(\omega)$
respectively; $v_n$ and $v_\lambda$ are thus functions from $\Omega$ into
$\mathbb{R}$.

Let $\mathcal{F}=\{f:\Omega\longrightarrow\mathbb{R}\}$; the Shapley
operator $J$ from $\mathcal{F}$ to itself is then defined by
$f\rightarrow J(f)$, where $J(f)$ is the function from $\Omega$ to
$\mathbb{R}$ satisfying

\begin{eqnarray}\label{Shapley}
J(f)(\omega)&=&\max_{\sigma\in\Delta(U)}\min_{\tau\in\Delta(V)}\left\{g(\sigma,\tau,\omega)+\sum_{\omega'\in\Omega}
f(w')\rho(\omega'|\sigma,\tau,\omega)\right\}\\
&=&\min_{\tau\in\Delta(V)}\max_{\sigma\in\Delta(U)}\left\{g(\sigma,\tau,\omega)+\sum_{\omega'\in\Omega}
f(w')\rho(\omega'|\sigma,\tau,\omega)\right\}
\end{eqnarray}

Then $J$ is nonexpansive on $\mathcal{F}$ endowed with the uniform
norm. The value $v_n$ of the $n$-stage game (resp. the value
$v_\lambda$ of the $\lambda$-discounted game) satisfies relation
(\ref{defvn}) (resp. (\ref{defvlambda})).

This recursive structure holds in a wide class of zero-sum repeated games and the
study of the asymptotic behavior of $v_n$ (resp. $v_\lambda$) as $n$
tends to $+\infty$ (resp. as $\lambda$ tends to 0) is a major topic
in game theory (see \cite{So2}
 for example)%: one wonders how values of a repeated game behave when it is played
%during a long time, either because of a high number of stages or
%because of a low discount factor
. Convergence of both $v_n$ and $v_\lambda$ (as well as equality of
the limits) has been obtained for different class of games, for
example absorbing games \cite{Koh}, recursive games \cite{Ev}, games
with incomplete information \cite{AuMa}, finite stochastic games
\cite{BeKo} \cite{BeKo2}, and Markov Chain Games with incomplete
information\cite{Ren}.

Even in the simple case of a finite stochastic game where the space
$\mathcal{F}$ on which $J$ is defined is $\mathbb{R}^n$, the Shapley
operator $J$ is only nonexpansive for the uniform norm
$\ell^\infty$. In the case of a general Shapley operator $J$, the
Banach space (which may be infinite dimensional) on which $J$ is
nonexpansive is always a set of bounded real functions (defined on a
set $\Omega$ of states) endowed with the uniform norm. As shown in
\cite{GuKe} and \cite{KoNe}, this lack of geometrical smoothness
implies that the families $v_n$ and $v_\lambda$ may not converge.
They may also converge to two different limits\cite{LeSo}. However
the goal of the so called "Operator Approach" (see \cite{RoSo} and
\cite{So}) is to infer, from specific properties in the framework of
games, convergence of both $v_n$ and $v_\lambda$ as well as equality
of their limits.

A closely related application, in the framework of discrete event
systems, is the problem of existence of the cycle-time of a topical
mapping \cite{GaGu} \cite{Gu}.

\subsection{Associated evolution equations}

In the current paper we investigate a slightly different direction :
the aim is to show that the sequence $v_n$ and the family
$v_\lambda$ defined in equations (\ref{defvn}) and
(\ref{defvlambda}) behave asymptotically as the solutions of certain
continuous-time evolution equations. This is interesting for at
least three reasons: first, this implies that proving the
convergence of $v_n$ or $v_\lambda$ reduces to study the asymptotic
of the solution of some evolution equation. Second, even if the
definitions (\ref{defvn}) of $v_n$ and (\ref{defvlambda}) of
$v_\lambda$ may seem dissimilar since one is recursive and the other
is a fixed point equation, we will see that the corresponding
equations in continuous time are of the same kind, hence it gives an
insight on the equality $\lim v_n=\lim v_\lambda$, satisfied for a
wide class of games. Third, we will prove in the process some
results of interest in their own right.

%Since we are primarily
%interested in the special case of a Shapley operator, which is nonexpansive for the uniform norm, we will never make any
%geometrical assumption on the underlying Banach space.
%For the same
%reason we will never assume anything on the mapping $J$ that is not
%satisfied in the case of Shapley operators.

\noindent Notice that equation (\ref{defVn}) can also be written as
a difference equation

\begin{equation}
(V_{n+1}-V_n)=J(V_n)-V_n
\end{equation}
\noindent which can be viewed as a discrete version of the evolution
equation

\begin{equation}\label{equadiff}
U'(t)=J(U(t))-U(t).
\end{equation}

\noindent Similarly, equations (\ref{defvn}) and
(\ref{eqrecvlambda}) can be considered as discrete versions of

\begin{equation}\label{equdiff2}
u'(t)=\Phi\left(\frac{1}{t+1},u(t)\right)-u(t)
\end{equation}
\noindent and
\begin{equation}\label{equdiff3}
u'(t)=\Phi(\lambda,u(t))-u(t)
\end{equation}
\noindent respectively. Notice that while (\ref{equdiff3}) is
autonomous, (\ref{equdiff2}) is not.

 The asymptotic relation between solutions of
(\ref{equadiff}) and (\ref{defVn}) will be discussed in section 3.
In that section we will a also prove some results about Eulerian
schemes related to (\ref{equadiff}), which have an interpretation in
terms of games with uncertain duration \cite{Ne,NeSo} in the case of
a Shapley Operator.

In section 4, we will study the asymptotic behavior of solutions of
the non-autonomous evolution equation

\begin{equation}\label{eqdiffPhi}
u'(t)=\Phi(\bm{\lambda}(t),u(t))-u(t)
\end{equation}
\noindent for some time-dependent parametrizations
$\bm{\lambda}(t)$, which in particular will cover both cases of
equations (\ref{equdiff2}) and (\ref{equdiff3}). We will first prove
that when $\bm{\lambda}(t)=\frac{1}{t}$  the solution of
(\ref{eqdiffPhi}) has the same asymptotic behavior,as $t$ goes to
$+\infty$, as the sequence $v_n$ as $n$ goes to $+\infty$. We will
then examine the case where the parametrization $\bm{\lambda}(t)$
converges slowly enough to 0, establishing that the solution of
(\ref{eqdiffPhi}) has then the same asymptotic behavior as the
family $v_\lambda$ as $\lambda$ goes to 0. Finally, using our
results in continuous time, we will study other dynamics in discrete
time generalizing (\ref{defvn}) and (\ref{eqrecvlambda}). Similarly
to section \ref{sectionVn}, in the case of a Shapley operator these
dynamics have an interpretation in terms of games with uncertain
duration.

\section{Dynamical system related to the operator $J$}\label{sectionVn}

%We recall that $(X,\|\cdot\|)$ is a Banach space and $J$ a
%non-expansive operator from $X$ to itself.
Let us denote $A=I-J$;
%, and for $\lambda\in\mathbb{R}$,
%$J_\lambda=I-\lambda A=\lambda J+(1-\lambda)I$.
 the operator $A$ is $m$-accretive, meaning that for any $\lambda>0$
 both properties are satisfied:
\begin{enumerate}
\item[$\mathrm{(i)}$] $\|x-y+\lambda A(x)-\lambda A(y)\|\geq\|x-y\|$ for all $(x,y)\in X^2$.
\item[$\mathrm{(ii)}$] $I+\lambda A$ is surjective.
\end{enumerate}

This implies that $A$ is maximal monotone\cite{Ka}.
% and $J_\lambda$
%can be seen as an Eulerian resolvent (compare to the usual proximal
%resolvent $(I+\lambda A)^{-1}$, see \cite{Br} for example). Note
%that for $\lambda\in[0,1]$ $J_\lambda$ is non-expansive, but not
%strictly contracting in general.
%We consider the family $V_n$ defined in (\ref{defVn}).
Recall that the analogous in continuous time of equation
(\ref{defVn}) defining $V_n$ is evolution equation (\ref{equadiff}),
%\begin{equation} \label{equadiff} U(t)+U'(t)=J(U(t))\end{equation}
which can also be written as

\begin{equation}\label{equdiffA}U'(t)=-A(U(t))\end{equation}

\noindent with initial condition $U(0)=U_0$, the Cauchy-Lipschitz
theorem ensuring the existence and uniqueness of such a solution.

\begin{example}\label{exempletrivialsuite}
Following example \ref{exempletrivial}, suppose $J(x)=x+c$. Then one
has $A(x)=-c$, so $U(t)=U_0+ct$.
\end{example}

%\noindent We are interested in the solution of an evolution equation
%involving a maximal monotone operator, but notice that there are two
%differences with the usual theory. First recall that no geometry
%asumption is made about the Banach space $X$. Also we cannot assume
%that $J$ has a fixed point, since this is usually not the case when
%$J$ is a Shapley operator, so no assumption is made about the
%non-emptiness of $A^{-1}(0)$.
This simple example shows that, as in discrete time where the true
sequence to consider is not $V_n$ but the normalized $v_n$, we are
not expecting convergence of $U(t)$ but rather of the normalized
quantity $\frac{U(t)}{t}$. This is a consequence of the fact that we
do not assume non emptiness of $A^{-1}(0)$.

Apart from equation (\ref{defVn}), there are numerous other natural
discretizations of equation (\ref{equdiffA}). For every $x_0 \in X$
and any sequence $\{\lambda_n\}$ in [0,1]\footnote{Usually these
schemes are defined for any sequence of positive steps, but here,
since we need the operators $I-\lambda_{n}A$ to be non expansive, we
have to assume that the $\lambda_n$ lie in $[0,1]$} the explicit
Eulerian scheme is defined by

\begin{equation}\label{eqeuler}
x_n-x_{n-1}=-\lambda_n A x_{n-1}\end{equation}that is

\begin{equation}
%x_n=J_{\lambda_n}x_{n-1}=\left(\prod_{i=n}^{1}
%[I-\lambda_{i}A]\right)(x_0)
x_n=\left(\prod_{i=n}^{1} [I-\lambda_{i}A]\right)(x_0).
\end{equation}
Notice that choosing $x_0=0$ and $\lambda_n=1$ for all $n$ leads to
the definition (\ref{defVn}) of $V_n$.

Other discrete trajectories are implicit proximal schemes (first
introduced when $A=\partial f$ in \cite{Mo}) which satisfy:
\[x_n-x_{n-1}=-\lambda_n A x_{n}\]
\noindent that is
\[x_n=\left(\prod_{i=n}^{1}
[I+\lambda_{i}A]^{-1}\right)(x_0).\]

\noindent In both cases we denote
\begin{eqnarray}
\label{defsigma}\sigma_n&=&\sum_{i=1}^n \lambda_i \\
\label{deftau}\tau_n&=&\sum_{i=1}^n \lambda^2_i.
\end{eqnarray}

Usually proximal schemes share better asymptotic properties (take
the simple example where $A$ is a rotation in $\mathbb{R}^2$ and
$\lambda_n\notin\ell^2$: then the proximal scheme will converge to
the fixed point of the rotation, while the Eulerian one will
diverge).
%In the case where $A=I-J$, we will show that some
%results involving proximal schemes (exponential formula and
%Kobayashi inequality) have an Eulerian counterpart.
However Eulerian schemes have the remarkable feature that they can
be computed explicitly, and they arise naturally in the
game-theoretic framework:

\begin{example}\label{uncertain}
When $J$ is the Shapley operator of a stochastic game $\Gamma$,
$x_n$ defined by (\ref{eqeuler}) is the non-normalized value of the
following $n-$stage game: states, actions, payoff and transition are
as in $\Gamma$, but at stage 1 there is a probability $1-\lambda_n$
that the game goes on to stage 2 without any payoff or transition.
Similarly at stage 2, there is no payoff nor transition with
probability $1-\lambda_{n-1}$, and at stage $n$ with probability
$1-\lambda_1$. In that case $\sigma_n$ and $\tau_n$ have a nice
interpretation: the expected number of stages really played is
$\sigma_n$, and the variance is $\sigma_n-\tau_n$. It is also
worthwile to notice that such games are particular cases of
stochastic games with uncertain duration\cite{Ne,NeSo}.
\end{example}

For this reason we will study exclusively Eulerian schemes, in the
case of an operator  $A=I-J$.
% they share well-known properties
%(exponential formula, Kobayashi inequality) of proximal
%trajectories.
Results of this section will be of three kind: first we study the
relative behavior of continuous and discrete dynamics when time goes
to infinity. Given a sequence $\lambda_n\notin\ell^1$ one
investigates the asymptotic relation between $U(\sigma_n)$ and the
n-th term $x_n$ of the Eulerian scheme defined in (\ref{eqeuler}).
This is done first in the special case of $V_n$ (Corollary
\ref{convvn}) and then in general (Corollary
\ref{generalasymptotic}).

We also consider the case of a fixed time $t$. In that case one cuts
the interval $[0,t]$ in a finite number $m$ of intervals of length
$\lambda_i$. These steps define an explicit scheme by
(\ref{eqeuler}), hence an approximate trajectory by linear
interpolation. One expects such a trajectory to be asymptotically
closer to the continuous trajectory defined by (\ref{equdiffA}) as
the discretization of the interval becomes finer. This is proved
first in the case where $\lambda_i=\frac{t}{m}$ for $1\leq i\leq m$
(Proposition \ref{expo}), and then generalized in Proposition
\ref{generalfixedtime}.

In the process we prove that two classical results, involving
proximal schemes and holding for any maximal monotone operator, have
an Eulerian counterpart when $A$ is of the form $I-J$: we establish
a new exponential formula in Proposition \ref{expo} and a
Kobayashi-like inequality in Proposition \ref{Kobaeuler}.
\subsection{Asymptotic study of the trajectory defined by equation
(\ref{equdiffA})}

The study of the asymptotic behavior of the solution of equation
(\ref{equadiff}) in general Banach spaces has started in the early
70's, in particular the main result of this subsection, Corollary
\ref{convvn} relating $v_n$ and $\frac{U(n)}{n}$, is already known
(see \cite{MiOh} and \cite{Ba}). Here we prove it in a different
way, similar to the first chapter of \cite{Br}, establishing during
the proof some inequalities that will be helpful in the remaining of
the paper.

\noindent Let us begin by proving several useful lemmas:

%\begin{lem}\label{udecr} Si une fonction $U$ $\mathcal{C}^1$ vérifie $\|U(t)+\mu U'(t)\|\leq\|U(t)\|$ pour un $\mu>0$, alors $\|U(t)\|$ est décroissante.
%\end{lem}

%\begin{proof}
%$U(t)+h U'(t)=\frac{h}{\mu}(U(t)+\mu u'(t))+(1-\frac{h}{\mu})u(t)$ et on a donc
%$\|u(t)+h u'(t)\|\leq\|u(t)\|$ $\forall h\leq\mu$. \\
%Soit $\varepsilon>0$, soit $f_\varepsilon(t)=\|u(t)\|-\varepsilon t$ et soit $a=\sup\{T\in\mathbb{R}^+,f_\varepsilon$ décroissante sur $[0,T]$\} et supposons $a$ fini. $u'$ étant uniformément continue sur $[a,a+\mu]$, $\exists\ 0<\eta<\mu$, $\forall (x,y)\in[a,a+\eta]^2,\|u'(x)-u'(y)\|\leq\varepsilon$. Pour $a\leq x< y\leq a+\eta$, on a alors
%\begin{eqnarray*}
%f_\varepsilon (y)-f_\varepsilon (x)&=&\|u(x)+(y-x) u'(z)\|-\|u(x)\|-\varepsilon (y-x) \textrm{ avec $z\in]x,y[$}\\
%&\leq &\|u(x)+(y-x) u'(x)\|+\varepsilon (y-x)-\|u(x)\|-\varepsilon (y-x)\\
%&\leq&0
%\end{eqnarray*}
%ce qui contredit la définition de a. On en déduit donc que $\forall \varepsilon>0, \|u(t)\|-\varepsilon t$ est décroissante sur $\mathbb{R}^+$, et donc que $\|u(t)\|$ l'est également.
%\end{proof}

\begin{lem}\label{Gronpart}
Let $f$ be a continuous function from $[a,b]\subset\mathbb{R}$ to
$\mathbb{R}$ such that for every $t\in[a,b]$

\[
f(t)\leq M+\int_a^t [g(s)+\beta(s)f(s)] ds
\]
\noindent for some continuous function $g$ and some non-negative
measurable function $\beta$ such that $\int_a^b \beta(s)ds<+\infty$.\\
Then $f$ satisfies
\[
f(t)\leq e^{\int_a^t \beta(s) ds}\left(M+\int_a^t g(s) e^{-\int_a^s
\beta(r) dr} ds\right)
\]
for all $t\in[a,b]$.
\end{lem}

\begin{proof}
Define $\alpha(t)=M+\int_a^t g(s)ds$. Since $f(t)\leq
\alpha(t)+\int_a^t \beta(s)f(s) ds$, Gronwall's inequality(\cite{Wa}
p. 15) implies that
\[
f(t)\leq \alpha(t)+ e^{\int_a^t \beta(s) ds }\cdot \int_a^t
\alpha(s) \beta(s) e^{-\int_a^s \beta(r) dr}ds.
\]
Integrating by part the last integral gives

\begin{eqnarray*}
f(t)&\leq& \alpha(t)+ e^{\int_a^t \beta(s) ds
}\cdot\left(\alpha(0)-\alpha(t)e^{-\int_a^t \beta(s) ds }+ \int_a^t
\alpha'(s) e^{-\int_a^s \beta(r) dr}ds\right)\\
&=&e^{\int_a^t \beta(s) ds}\left(M+\int_a^t g(s) e^{-\int_a^s
\beta(r) dr} ds\right).
\end{eqnarray*}
\end{proof}

In the remaining of the paper we will repeatedly use the following
consequence of Lemma \ref{Gronpart}:
\begin{proposition}\label{decrexpo}
If $y:[a,b]\subset\mathbb{R}\rightarrow X $ is an absolutely
continuous function satisfying for every $t\in[a,b]$

\[\|y(t)+y'(t)\|\leq(1-\gamma(t))\|y(t)\|+h(t)\]

\noindent where $\gamma$ is a continuous function from $[a,b]$ to
$[-\infty,1]$ and $h$ is a continuous function from $[a,b]$ to
$\mathbb{R}$, then $y$ satisfies

\[
\|y(t)\|\leq\ e^{-\int_{a}^{t}\gamma(s) ds}\left(\|y(a)\|+\int_a^t
h(s)e^{\int_a^s \gamma(r) dr} ds\right)
\]

\noindent for all $t\in[a,b]$.
\end{proposition}

\begin{proof}
$z(t)=y(t)e^t$ satisfies
$\|z'(t)\|\leq(1-\gamma(t))\|z(t)\|+h(t)e^t$, hence for every
$t\in[a,b]$

\begin{eqnarray*}\|z(t)\|&\leq& \|z(a)\|+\left\|\int_a^t z'(s)
ds\right\|\\&\leq& \|z(a)\|+\int_a^t \|z'(s)\| ds \\
&\leq& \|z(a)\|+\int_a^t h(s)e^s+(1-\gamma(s))
\|z(s)\|.\end{eqnarray*} Applying Lemma \ref{Gronpart} to $\|z\|$
thus gives
\[
\|z(t)\|\leq e^{\int_a^t 1-\gamma(s) ds}\left(\|z(a)\|+\int_a^t
h(s)e^s e^{-\int_a^s 1-\gamma(r) dr} ds\right).
\]
Multiplying each side by $e^{-t}$ implies the result.
%Applying Gronwall's lemma to the function $\|r\|$ we thus deduce
%that \[\|r(y)\|\leq\|r(x)\|\cdot e^{\int_{x}^{y}1-\gamma(t) dt}.\]
%Multiplying each side of this inequality by $e^{-y}$ we get the
%result.
\end{proof}

We now use this technical result to compare two solutions of
(\ref{equadiff}):
\begin{proposition} If both $U$ and $V$ satisfy (\ref{equadiff}), then $\|U(t)-V(t)\|$ is non-increasing.
\end{proposition}

\begin{proof} Define $f=U-V$ which satisfies

\begin{equation*}
\|f(t)+f'(t)\|= \|J(U(t))-J(V(t))\|\leq\|U(t)-V(t)\|=\|f(t)\|.
\end{equation*}
Apply the preceding proposition to $\gamma\equiv0$ and $f$.
\end{proof}

\begin{corollary}\label{decr} If $U$ is a solution of (\ref{equadiff}), then $\|U'(t)\|$ is non-increasing.
\end{corollary}

\begin{proof} Let $h>0$ and $U_h(t)=U(t+h)$. The function $U_h$ satisfies equation (\ref{equadiff}), so applying the preceding proposition to $U$ and $U_h$ we get that $t\rightarrow\frac{\|U(t+h)-U(t)\|}{h}$ is non-increasing on $\mathbb{R}^+$.
Letting $h$ go to 0 gives the result
\end{proof}

%We now prove that the following lemma, initially obtained in the
%Hilbert framework, holds also here. The crucial argument, as in the
%Hilbert case, is that the norm of the derivative of the solution of
%(\ref{equadiff2}) is non-increasing, and the proof is similar.
An interesting consequence of Corollary \ref{decr} is the following
inequality, proved in Chapter 1 of \cite{Br}:
\begin{lem}[Chernoff's estimate] \label{Chernoff} Let $U$ be the solution of (\ref{equadiff}) with $U(0)=U_0$.
Then \[ \|U(t)-J^n(U_0)\|\leq\|U'(0)\|\sqrt{t+[n-t]^2}.\]
\end{lem}

\begin{proof}[Sketch of proof]
%It is enough to handle the case $\mu=1$ because of Remark
%\ref{changmtvar}.
Proceed by induction on $n$; the proof for the case $n=0$ comes from
the fact that $\|U'\|$ is non-increasing by Corollary \ref{decr}.

\end{proof}

In particular if we take $U_0=0$ and $t=n$ in Lemma \ref{Chernoff},
we finally get the following corollary relating continuous and
discrete trajectories:

\begin{corollary} \label{convvn}
The solution $U$ of (\ref{equadiff}) with $U(0)=0$ satisfies
\[\left\|\frac{U(n)}{n}-v_n\right\|\leq\frac{\|J(0)\|}{\sqrt{n}}\]
In particular $v_n$ converges iff $\frac{U(t)}{t}$ converges, and
then the limits are the same.

\end{corollary}

\begin{proof}
The only point that remains to be shown is that if the sequence
$\frac{U(n)}{n}$ converges as $n$ tends to $+\infty$, then so does
$\frac{U(t)}{t}$ as $t$ tends to $+\infty$.

Using Corollary \ref{decr}, we obtain
\[
\|U(t)-U([t])\|\leq(t-[t])\|U'(0)\|\leq \|U'(0)\|
\]
which implies that $\frac{U(t)}{t}-\frac{U([t])}{[t]}$ goes to 0 as
$t$ tends to $+\infty$.

\end{proof}
\subsection{An exponential formula}
%Two well-studied discretisations of equation (\ref{equadiff}) are
%the Eulerian (explicit) and proximal (implicit) schemes, the former
%defined by
%\[
% x_{n+1}\in(I-\lambda_n A)(x_n)
%\]
%and the latter, first introduced when $A=\partial f$ in \cite{Mo},
%by
%\[
%x_{n+1}\in(I+\lambda_n A)^{-1}(x_n)
%\]
%\noindent with $\lambda_n$ a sequence in $\mathbb{R}^+$.
%
%In this subsection and the following one we prove that some
%classical results about proximal schemes have an analoguous for
%Eulerian schemes when the maximal monotone operator $A$ equals $I-J$
%for some nonexpansive $J$.

When $A$ is a $m$-accretive operator on a Banach space, a
fundamental result (see \cite{CrLi} p. 267) is that the solution $U$
of (\ref{equadiff}) satisfies the following exponential formula for
every $t\geq0$, where the convergence is strong:
\begin{equation}\label{expprox}
\lim_{m\rightarrow +
\infty}\left(I+\frac{t}{m}A\right)^{-m}(U_0)=U(t)
\end{equation}
\noindent

In the special case where $J$ is a nonexpansive operator and
$A=I-J$, we now establish an Eulerian analogous of this classical
"proximal exponential formula".

\begin{déf}
For $x\in X$, $l\in\mathbb{N}$ and $t\in\mathbb{R}^+$, let us denote
\begin{equation}
U_t^m(x)=\left(I-\frac{t}{m}A\right)^m(x)
\end{equation}
the $m$-th term of an Eulerian scheme with steps $\frac{t}{m}$.

\end{déf}

\begin{proposition} \label{expo} Let $U_0\in X$ and $U$ the solution of (\ref{equadiff}) with $U(0)=U_0$. Then if $m\geq t$,
\begin{equation}\label{expeuler}
\|U_{t}^m(U_0)-U(t)\|\leq
\|A(U_0)\|\frac{t}{\sqrt{m}}.\end{equation} In particular, for any
$t\geq 0$, the following strong convergence holds:

\begin{equation}
\lim_{m\rightarrow +
\infty}\left(I-\frac{t}{m}A\right)^{m}(U_0)=U(t)
\end{equation}
\end{proposition}

%\begin{proof}
%For $\lambda\leq 1$ we apply Lemma \ref{Chernoff} to $J_\lambda$,
%noticing that $A=\frac{1}{\lambda}[I-J_\lambda]$ which gives:
%\[\|U(t)-(I-\lambda A)^n U_0)\|\leq\|U'(0)\|\sqrt{\lambda t+[n\lambda-t]^2}\]
%If $n\geq t$, we take $\lambda=\frac{t}{n}\leq 1$ and we obtain:
%\[\|U(t)-U_{t}^n(U_0)\|\leq \|A(U_0)\|\frac{t}{\sqrt{n}}\]
%\end{proof}

\begin{proof}
For any $\lambda\in[0,1]$, $J_\lambda:=\lambda J+(1-\lambda)
I=I-\lambda A$ is nonexpansive. Denote by $U_\lambda$ the solution
of
\begin{equation}\label{equadiff2}
U_\lambda(t)+U_\lambda'(t)=J_\lambda(U_\lambda(t))
\end{equation}
 \noindent with $U_\lambda(0)=U_0$. Applying Lemma \ref{Chernoff} to $U_\lambda$ and the nonexpansive operator
 $J_\lambda$:
\[ \|U_\lambda(t)-J_\lambda^n(U_0)\|\leq\|U_\lambda'(0)\|\sqrt{t+[n-t]^2}\]
\noindent so in particular for $n=t$
\begin{equation} \label{Chernofflambda}
\|U_\lambda(n)-J_\lambda^n(U_0)\|\leq\|U_\lambda'(0)\|\sqrt{n}.\end{equation}
Denote by $U$ the the solution of (\ref{equadiff}) with $U(0)=U_0$
and notice that the function $t\rightarrow U(\lambda t)$ satisfies
(\ref{equadiff2}) and has the same initial condition as $U_\lambda$.
This implies that $U_\lambda(t)=U(\lambda t)$ and putting this in
(\ref{Chernofflambda}),
\[ \|U(\lambda n)-J_\lambda^n(U_0)\|\leq\lambda\|U'(0)\|\sqrt{n}.\]
For any $t'\leq n$, choosing $\lambda=\frac{t'}{n}\in[0,1]$ thus
gives
\[\left\|U(t')-\left(I-\frac{t'}{n}A\right)^{n}(U_0)\right\|\leq\|U'(0)\|\frac{t'}{\sqrt{n}}.\]
\noindent which is the desired result.
\end{proof}

%\begin{remark} The formula (\ref{expeuler}) can be seen as an Eulerian
%analogous of the classical exponential formula
%\begin{equation}\label{expprox}
%\lim_{m\rightarrow +
%\infty}\left(I+\frac{t}{m}A\right)^{-m}(U_0)=U(t)
%\end{equation}
%\noindent (see \cite{CrLi} for example). Notice however that formula
%(\ref{expprox}) holds for any maximal monotone operator, while
%proposition \ref{expo} relies on the nonexpansiveness of $J$.
%\end{remark}
\subsection{Comparaison of two Eulerian schemes}
%In the preceding section we proved that a classical proximal result
%which holds for maximal monotone operators in Banach spaces has an
%explicit Eulerian counterpart in the case where $A=I-J$. In the same
%vein we prove that the Kobayashi inequality for two proximal schemes
%established in \cite{Ko} also holds for explicit schemes in this
%special case.
%\begin{remark}\label{remarquetoutebete}
%If $\{\lambda_n\}$ is the constant sequence equal to $\frac{t}{m}$,
%then $x_m=U_t^m(x_0)$.\\
%If $\{\lambda_n\}$ is the constant sequence equal to $1$ and
%$x_0=0$, then $x_m=V_m$.
%\end{remark}
%The previous example indicates that the normalized quantity to
%consider is $\displaystyle{\frac{x_n}{\sigma_n}}$. The aim of this
%subsection and of the following ones will be to compare it to the
%normalized quantity in continuous time
%$\displaystyle{\frac{U(\sigma_n)}{\sigma_n}}$.
To generalize Proposition \ref{expo} to explicit schemes with
arbitrary steps, it is useful to estimate first the difference
between two Euler schemes: let $x_0$ and $\hat{x}_0$ in $X$,
$\{\lambda\}_n$ and $\{\hat{\lambda}_n\}$ two sequences in ]0,1].
Define $x_n$ , $\sigma_n$ and $\tau_n$ (resp. $\hat{x}_n$,
$\hat{\sigma}_n$ and $\hat{\tau_n}$) as in (\ref{eqeuler}),
(\ref{defsigma}) and (\ref{deftau}). The following proposition,
which gives a majoration of the distance between two Eulerian
trajectories, is an analogous of the classical Kobayashi inequality
(Lemma 2.1 in \cite{Ko}) which gives a majoration of the distance
between two proximal trajectories:
%, et $\Lambda$ et $\hat{\Lambda}$ les bornes
%des suites $\{\lambda_n\}$ et $\{\hat{\lambda}_n\}$ (on rappelle
%qu'on suppose que $\Lambda$ et $\hat{\Lambda}$ sont inférieurs à 1)

\begin{proposition}\label{Kobaeuler}
For any $z\in X$ and $(k,l)\in\mathbb{N}^2$,
\[\|x_k-\hat{x}_l\|\leq \|x_0-z\|+\|\hat{x}_0-z\|+\|A(z)\|\sqrt{(\sigma_k-\hat{\sigma}_l)^2+\tau_k+\hat{\tau}_l}\]
\end{proposition}

\begin{proof}
We proceed by induction and begin by the case $l=0$.\\
We recall that $J_\lambda=I_\lambda A$ is non-expansive for
$\lambda\leq 1$, so we obtain
\begin{eqnarray*}
\|x_j-z\|& \leq & \|x_j-J_{\lambda_j}(z)\|+\|J_{\lambda_j}(z)-z\| \\
 & = & \|J_{\lambda_j}(x_{j-1})-J_{\lambda_j}(z)\|+\lambda_j\|A(z)\|\\
 & \leq & \|x_{j-1}-z\|+\lambda_j\|A(z)\|
\end{eqnarray*}
and summing these inequalities for $i\in\{1,\cdots,k\}$ we get
$\|x_k-z\|\leq\|x_0-z\|+\sigma_k\|A(z)\|$, which implies that
\[\|x_k-\hat{x}_0\|\leq \|x_k-z\|+\|\hat{x}_0-z\|\leq\|x_0-z\|+\|\hat{x}_0-z\|+\sigma_k\|A(z)\|\]
and the proposition holds when $l=0$. The case $k=0$ is proved in the same way. \\
We will now assume the formula to be true for $(k-1,l)$, $(k,l-1)$
et $(k-1,l-1)$ and deduce that it also holds for $(k,l)$.

Define numbers
$\alpha_{k,l}=\frac{\lambda_k(1-\hat{\lambda}_l)}{\lambda_k+\hat{\lambda}_l-\lambda_k
\hat{\lambda}_l}$, $\beta_{k,l}=\frac{\hat{\lambda}_l
(1-{\lambda_k})}{\lambda_k+\hat{\lambda}_l -\lambda_k
\hat{\lambda}_l}$ et $\gamma_{k,l}=\frac{\lambda_k
\hat{\lambda}_l}{\lambda_k+\hat{\lambda}_l -\lambda_k
\hat{\lambda}_l}$ and note that they are non-negative with sum 1.
Introduce also
$c_{k,l}=\sqrt{(\sigma_k-\hat{\sigma}_l)^2+\tau_k+\hat{\tau}_l}$.
For any $x$ and $y$ in $X$, one check that the following equality
holds:
\[J_{\lambda_k}(x)-J_{\hat{\lambda}_l}(y)=\alpha_{k,l}(J_{\lambda_k}(x)-y)+\beta_{k,l}(x-J_{\hat{\lambda}_l}(y))+\gamma_{k,l}(J(x)-J(y)).\]
In particular, letting $x=x_{k-1}$, $y=\hat{x}_{l-1}$ and using the
non-expansiveness of $J$, we get
\[\|x_k-\hat{x}_l\|\leq \alpha_{k,l}\|x_k-\hat{x}_{l-1}\|+\beta_{k,l}\|x_{k-1}-\hat{x}_{l}\|+ \gamma_{k,l}\|x_{k-1}-\hat{x}_{l-1}\|\]
so by induction,
\begin{eqnarray*}
\|x_k-\hat{x}_l\|&\leq& \|x_0-z\|+\|\hat{x}_0-z\|+ \|A(z)\|(\alpha_{k,l}c_{k,l-1}+\beta_{k,l}c_{k-1,l}+\gamma_{k,l}c_{k-1,l-1})\\
&\leq&\|x_0-z\|+\|\hat{x}_0-z\|+ \|A(z)\|\sqrt{\alpha_{k,l}+\beta_{k,l}+\gamma_{k,l}} \sqrt{d_{k,l}}\\
&=& \|x_0-z\|+\|\hat{x}_0-z\|+ \|A(z)\| \sqrt{d_{k,l}}
\end{eqnarray*}
where we have denoted $d_{k,l}=\alpha_{k,l}c_{k,l-1}^2+\beta_{k,l}c_{k-1,l}^2+\gamma_{k,l}c_{k-1,l-1}^2$.\\
In addition,
\begin{eqnarray*}
c_{k,l-1}^2&=&(\sigma_k-\hat{\sigma}_{l-1})^2+\tau_k+\hat{\tau}_{l-1}\\
&=&(\sigma_k-\hat{\sigma}_{l}+\hat{\lambda}_l)^2+\tau_k+\hat{\tau}_{l-1}\\
&=&(\sigma_k-\hat{\sigma}_l)^2+\hat{\lambda}_l^2+2\hat{\lambda}_l(\sigma_k-\hat{\sigma}_{l})+ \tau_k+\hat{\tau}_{l-1}\\
&=&c_{k,l}^2+2\hat{\lambda}_l(\sigma_k-\hat{\sigma}_{l})
\end{eqnarray*}
and similarly,
\[ c_{k-1,l}^2=c_{k,l}^2-2\lambda_k(\sigma_k-\hat{\sigma}_{l}).\]
Moreover
\begin{eqnarray*}
c_{k-1,l-1}^2&=&(\sigma_{k-1}-\hat{\sigma}_{l-1})^2+\tau_{k-1}+\hat{\tau}_{l-1}\\
&=&(\sigma_k-\hat{\sigma}_{l}+\hat{\lambda}_l-\lambda_k)^2+\tau_{k-1}+\hat{\tau}_{l-1}\\
&=&(\sigma_k-\hat{\sigma}_l)^2+\hat{\lambda}_l^2+\lambda_k^2+2(\hat{\lambda}_l-\lambda_k)(\sigma_k-\hat{\sigma}_{l})-2\hat{\lambda}_l \lambda_k+\tau_{k-1}+\hat{\tau}_{l-1}\\
&=&c_{k,l}^2+2(\hat{\lambda}_l-\lambda_k)(\sigma_k-\hat{\sigma}_{l})-2\hat{\lambda}_l
\lambda_k.
\end{eqnarray*}
So
\begin{eqnarray*}
d_{k,l}&=&\alpha_{k,l}c_{k,l-1}^2+\beta_{k,l}c_{k-1,l}^2+\gamma_{k,l}c_{k-1,l-1}^2\\
&=&c_{k,l}^2+2(\sigma_k-\hat{\sigma}_{l})(\alpha_{k,l}\hat{\lambda}_l-\beta_{k,l}\lambda_k+\gamma_{k,l}(\hat{\lambda}_l-\lambda_k))-2\hat{\lambda}_l \lambda_k \gamma_{k,l}\\
&=&c_{k,l}^2+2\frac{\sigma_k-\hat{\sigma}_{l}}{\lambda_k+\hat{\lambda}_l-\lambda_k \hat{\lambda}_l}(\lambda_k\hat{\lambda}_l(1-\hat{\lambda}_l)-\lambda_k\hat{\lambda}_l(1-\lambda_k)+\lambda_k\hat{\lambda}_l(\hat{\lambda}_l-\lambda_k))-2\hat{\lambda}_l \lambda_k \gamma_{k,l}\\
&=&c_{k,l}^2-2\hat{\lambda}_l \lambda_k \gamma_{k,l}\\
&\leq&c_{k,l}^2
\end{eqnarray*}
and we have established that
\[\|x_k-\hat{x}_l\|\leq \|x_0-z\|+\|\hat{x}_0-z\|+\|A(z)\|\sqrt{(\sigma_k-\hat{\sigma}_l)^2+\tau_k+\hat{\tau}_l}.\]

\end{proof}

\subsection{Comparaison of an Eulerian scheme to a continuous
trajectory}

We now combine the results of the two preceding subsections:
Proposition \ref{expo} comparing the continuous trajectory with a
particular Eulerian scheme, and Proposition \ref{Kobaeuler} relating
any two Eulerian schemes.
% Using Remark \ref{remarquetoutebete},
%We deduce from Proposition \ref{Kobaeuler} that if $x_k$ is an
%Eulerian trajectory,

%By letting $m$ goes to $+\infty$ and using , we deduce the

\begin{corollary}\label{Cordiscrcont} Let $\{x_n\}_{n\in\mathbb{N}}$ be an Eulerian scheme as defined in (\ref{eqeuler}). Then for any $t\geq 0$ and $k\in\mathbb{N}$,

\[\|x_k-U(t)\|\leq \|x_0-U_0\|+\|A(U_0)\|\sqrt{(\sigma_k-t)^2+\tau_k}\]

\end{corollary}

\begin{proof}
Apply Proposition \ref{Kobaeuler} to $x_k$ and $U_t^m(U_0)$ to get
\[\|x_k-U_t^m(U_0)\|\leq \|x_0-U_0\|+\|A(U_0)\|\sqrt{(\sigma_k-t)^2+\tau_k+\frac{t^2}{m}}.\]
Let $m$ go to $+\infty$ and use Proposition \ref{expo}.
\end{proof}

This corollary has some interesting consequences in two directions,
as it generalizes both Corollary \ref{convvn} and Proposition
\ref{expo}. First, it shows that any normalized discrete trajectory
behave as the normalized continuous one as time goes to infinity:

\begin{corollary}\label{generalasymptotic}
For any $t\geq 0$ and any Eulerian scheme $\{x_i\}$ such that
$\sigma_k=t$,
\begin{eqnarray*}\frac{\|x_k-U(t)\|}{t}&\leq&
\frac{\|x_0-U_0\|+\|A(U_0)\|\sqrt{t}}{t}
%&\leq&\frac{\|x_0-U_0\|}{\sigma_k}+\frac{\|A(U_0)\|}{\sqrt{\sigma_k}}\\
\end{eqnarray*}
\end{corollary}
\begin{proof}
Apply Corollary \ref{Cordiscrcont} and use the fact that
$\tau_k\leq\sigma_k$ since all $\lambda_i$ are in $[0,1]$.
\end{proof}
\vspace{.2cm}

On the other hand, take now the case of a fixed time $t$. Let $U$ be
the solution of (\ref{equdiffA}) with initial condition $U(0)=U_0$,
and let $\{x_i\}_{0\leq i\leq n}$ defined by (\ref{eqeuler}) be an
Eulerian scheme with same initial condition $x_0=U_0$ and
$\sigma_n=t$. One constructs an approximation $\widetilde{x}$ of the
continuous trajectory $U$ on the interval $[0,t]$ by
$\widetilde{x}(\sigma_k)=x_k$ for $0\leq k\leq n$, and linear
interpolation on intervals $[\sigma_k,\sigma_{k+1}]$. The following
proposition states that such approximation $\widetilde{x}$ will
becomes asymptotically close to $U$ as the discretization $0\leq
\lambda_1\leq \lambda_1+\lambda_2\leq\cdots\leq
\sigma_{n-1}\leq\sigma_n=t$ of the interval $[0,t]$ gets finer:

\begin{proposition}\label{generalfixedtime}
For any $t'$ in the interval $[0,t]$,
\[
\|\widetilde{x}(t')-U(t')\|\leq
\|A(U_0)\|(1+(1+\sqrt{2})t)\cdot\sqrt{\max_{1\leq i\leq
n}\{\lambda_i\}}.
\]
\end{proposition}

\begin{proof}
Let $t'\in[0,t]$ and $k$ such that $\sigma_{k-1}\leq t'\leq
\sigma_k$. Noticing that
\[
\|\widetilde{x}(t')-U(t')\|\leq
\|\widetilde{x}(t')-\widetilde{x}(\sigma_k)\|+\|\widetilde{x}(\sigma_k)-U(\sigma_k)\|+\|U(\sigma_k)-U(t')\|\\
\]
\noindent we will evaluate the three components of the right-hand
side separately.

Since $\widetilde{x}$ is affine on $[\sigma_k,\sigma_{k+1}]$,
applying Proposition \ref{Kobaeuler} gives
\begin{eqnarray}
\nonumber \|\widetilde{x}(t')-\widetilde{x}(\sigma_k)\|&\leq&\|\widetilde{x}(\sigma_{k-1})-\widetilde{x}(\sigma_{k})\|\\
\nonumber&=&\|x_k-x_{k-1}\|\\
\nonumber&\leq&\|A(U_0)\|\sqrt{(\sigma_k-\sigma_{k-1})^2+\tau_k+\tau_{k-1}}\\
\label{ineg1} &=&\|A(U_0)\|\sqrt{2\tau_k}.
\end{eqnarray}

On another hand, Corollary \ref{Cordiscrcont} implies that
\begin{eqnarray}
\nonumber\|\widetilde{x}(\sigma_k)-U(\sigma_k)\|&=&\|x_k-U(\sigma_k)\|\\
\label{ineg2} &\leq& \|A(U_0)\|\sqrt{\tau_k}.
\end{eqnarray}

Thirdly, using the mean value Theorem as well as Corollary
\ref{decr},
\begin{eqnarray}
\nonumber\|U(\sigma_k)-U(t')\|&\leq &|\sigma_k-t'|
\max_{t''\in[t',\sigma_k]}
\|U'(t'')\|\\
\nonumber&\leq&|\sigma_k-\sigma_{k-1}|\cdot\|U'(0)\|\\
\label{ineg3} &=&\lambda_k\|A(U_0)\| .
\end{eqnarray}
Adding inequalities (\ref{ineg1}) to (\ref{ineg3}) we thus deduce
that
\[
\|\widetilde{x}(t')-U(t')\|\leq \|A(U_0)\| (
\lambda_k+(1+\sqrt{2})\sqrt{\tau_k}).
\]
We use the facts that $\lambda_k\leq \sqrt{\lambda_k}\leq
\sqrt{\max_{1\leq i\leq n}\{\lambda_i\}}$, and that
$\tau_k\leq\tau_n\leq\ t \max_{1\leq i\leq n}\{\lambda_i\}$ to
conclude the proof.
\end{proof}

This proposition has an interpretation in the particular framework
of Example \ref{uncertain}: consider a game with an expected
duration of $t$. The previous result establishes that this game has
a non normalized value close to $U(t)$, providing that at each stage
the probability of playing is small (that is to say, if there is a
high variance in the number of stages really played).

\section{Dynamical systems linked to the family $\Phi(\lambda,\cdot)$}
%For $0<\lambda\leq1$, we denote $\Phi(\lambda,x)=\lambda J(\frac{1-\lambda}{\lambda}x)$. $\Phi(\lambda,\cdot)$ is $1-\lambda$ contracting, let $v_\lambda$ be its unique fixed point.\\
Let $\bm{\lambda}:\mathbf{R}\rightarrow ]0,1]$ be a continuous
function. In this section we study the asymptotic behavior of the
solution to evolution equation (\ref{eqdiffPhi}):
\begin{equation*} u(t)+u'(t)=\Phi(\bm{\lambda}(t),u(t)) \quad \mathrm{with}\
u(0)=u_0\end{equation*} where $\Phi$ is the operator defined by
equation (\ref{defPhi}).% To ease the reading, from now on
%$\bm{\lambda}$ will be written in bold when it is a function,
%whereas an element of ]0,1] will be denoted as $\lambda$.

\begin{remark}
Since the mapping $(x,t)\rightarrow \Phi(\bm{\lambda}(t),x)-x$ is
globally 2-Lipschitz in its first variable, Cauchy-Lipschitz-Picard
theorem ensures the existence and uniqueness of the solution of
(\ref{eqdiffPhi}), and that it is defined on the whole set
$\mathbb{R}^+$.
\end{remark}

%\begin{remark}
When the recession function $\Phi(0,\cdot)$ exists, any accumulation
point $v$ of $v_n$ or $v_\lambda$ will satisfy
\begin{equation}\label{Phi0v}
\Phi(0,v)=v
\end{equation}
but equation (\ref{Phi0v}) may have many solutions (for example in
the case of games with incomplete information \cite{RoSo} any
convex/concave function satisfies (\ref{Phi0v})). The evolution
equation (\ref{eqdiffPhi}) may thus be seen as a perturbation of
(\ref{Phi0v}), and we will study the effect of some perturbations on
the asymptotic behavior of the solution of (\ref{eqdiffPhi}). See
for example \cite{AtCo} for a similar approach in the framework of
convex minimization.

The main results of this section are the following:
\begin{itemize}
\item When $\bm{\lambda}$ is the constant $\lambda$, the solution of
(\ref{eqdiffPhi}) converges to $v_\lambda$.
\item When $\bm{\lambda}(t)\sim\frac{1}{t}$, the solution of
(\ref{eqdiffPhi}) behave asymptotically as the family $\{v_n\}$
\item When $\bm{\lambda}(t)$ converges to 0 slowly enough, the solution of
(\ref{eqdiffPhi}) behave asymptotically as the family
$\{v_\lambda\}$
%\end{remark}
\end{itemize}
\noindent The first two results are not surprising since in those
cases evolution equation (\ref{eqdiffPhi}) is a continuous version
of equation (\ref{eqrecvlambda}) or (\ref{defvn}) respectively. The
third result is of a different nature but is also natural: denote by
$u_\lambda$ the solution of (\ref{eqdiffPhi}) when $\bm{\lambda}$ is
the constant $\lambda$. We establish that if the parametrization
$\bm{\lambda}$ in (\ref{eqdiffPhi}) is of slow variation, the
solution $u$ evaluated at time $t$ is close to $u_{\lambda(t)}(t)$,
hence to $v_{\lambda(t)}$ (see figure below).

%\begin{center}
%\begin{figure}[H]\centering
%\input{dessinguillaume.pstex_t}
%%\caption{}
%\end{figure}
%\end{center}

In the process of proving those three results, we also answer
natural questions about the behavior of the solution $u$ of equation
(\ref{eqdiffPhi}) as a function of the parameters, namely we will
prove that:

\begin{itemize}
\item If $\bm{\lambda}\notin\ell^1$ the asymptotic behaviour of $u$ does not depend on the
initial value $u_0$.
\item If two parametrizations $\bm{\lambda}$ and $\bm{\widetilde{\lambda}}$
are asymptotically close, then it is also the case for the
corresponding solutions $u$ and $\widetilde{u}$.

\end{itemize}

First we prove a simple fact that will be repeatedly used in the
remaining of the paper. Recall, by equation (\ref{defvlambda}), that
for any $t\geq0$, $v_{\bm{\lambda}(t)}$ is the only solution of
\begin{equation}\label{defvlambdat}
v_{\bm{\lambda}
(t)}=\Phi\left(\bm{\lambda}(t),v_{\bm{\lambda}(t)}\right).
\end{equation}
The following Lemma relates the behavior of $u'(t)$ to that of
$u(t)-v_{\bm{\lambda}(t)}$:
\begin{lem}\label{simple}
Let $u$ be the solution of evolution equation (\ref{eqdiffPhi}) and
$v_{\bm{\lambda}(\cdot)}$ be defined by (\ref{defvlambdat}). Then
for any $t\geq0$, $\|u(t)-v_{\bm{\lambda}
(t)}\|\leq\frac{\|u'(t)\|}{\bm{\lambda} (t)}$
\end{lem}
\begin{proof}
\begin{eqnarray*}
\|u'(t)\|&=&\left\|u(t)-\Phi(\bm{\lambda}(t),u(t))\right\| \\
&\geq&\left\|u(t)-v_{\bm{\lambda} (t)}\right\|-\left\|\Phi(\bm{\lambda}(t),u(t))-\Phi(\bm{\lambda}(t),v_{\bm{\lambda} (t)})\right\| \\
&\geq&\left\|u(t)-v_{\bm{\lambda}
(t)}\right\|-\left(1-\bm{\lambda}(t)\right)\left\|u(t)-v_{\bm{\lambda} (t)}\right\|\\
 &=&\bm{\lambda} (t) \left\|u(t)-v_{\bm{\lambda} (t)}\right\|.
\end{eqnarray*}
\end{proof}

\subsection{Constant case}
We start by considering the simplest case where the function
$\bm{\lambda}$ is a constant $\lambda$. Equation (\ref{eqdiffPhi})
is then a continuous analogous of equation (\ref{eqrecvlambda}), so
one can expect that $u(t)$ converges to $v_\lambda$, and indeed this
is the case.

Start by a technical lemma:

\begin{lem}\label{simple2}
If f satisfies $f(t)+f'(t)=B(f(t))$, where $B$ is an $1-\lambda$
contracting operator, then

\[\|f'(t)\|\leq \|f'(0)\|\cdot e^{-\lambda t}.\]
\end{lem}

\begin{proof}
Let $h>0$ and $\displaystyle{f_h(t)=\frac{f(t+h)-f(t)}{h}}$. Since
$B$ is $(1-\lambda)$ contracting:
\begin{eqnarray}
\|f_h(t)+f_h'(t)\|&=&\frac{1}{h}\|f(t+h)+f'(t+h)-[f(t)+f'(t)]\|\\
&=&\frac{1}{h}\|B(f(t+h))-B(f(t))\|\\
&\leq&(1-\lambda)\|f_h(t)\|.
\end{eqnarray}

Proposition \ref{decrexpo} applied to $f_h$ thus implies that
\[\|f_h(t)\|\leq \|f_h(0)\|\cdot e^{-\lambda t}\]
and letting $h$ go to 0 gives the result.

\end{proof}

An immediate consequence is:
\begin{corollary}\label{lambdaconstant}
If $u$ is the solution of (\ref{eqdiffPhi}) with
$\bm{\lambda}(t):=\lambda$, then
\[\lim_{t\rightarrow+\infty}u(t)=v_\lambda\]
\end{corollary}
\begin{proof}
Lemmas \ref{simple} and \ref{simple2} imply that
\[\|u(t)-v_{\lambda}\|=\|u(t)-v_{\bm{\lambda}
(t)}\|\leq\frac{\|u'(t)\|}{\bm{\lambda}
(t)}\leq\|u'(0)\|\cdot\frac{e^{-\lambda t}}{\lambda}\] and the right
member goes to 0 as $t$ tends to $+\infty$.
\end{proof}
\subsection{Some generalities on the non-autonomous case}
%\textrm{}
%\\

The case when the parametrization $\bm{\lambda}$ is not constant is
more difficult to handle: the same method as in the proof of
corollary \ref{lambdaconstant} leads to
\[u(t+h)-u(t)+u'(t+h)-u'(t)=\Phi(\bm{\lambda}(t+h),u(t+h))-\Phi(\bm{\lambda}(t),u(t))\]
but Proposition \ref{decrexpo} does not apply.

However, we can prove if the perturbation is strong enough:
\begin{proposition}\label{depart}
If $\int_{0}^{+\infty}\bm{\lambda}(t) dt=+\infty$, the asymptotic
behavior of $u$ solution of (\ref{eqdiffPhi}) does not depend of the
choice of $u(0)$.
\end{proposition}
\begin{proof}
Let $u$ and $v$ be two solutions of (\ref{eqdiffPhi}), define the
function $g$ by $g(x)=\|u(x)-v(x)\|$. According to proposition
\ref{decrexpo},
\[g(x)\leq g(0)\cdot e^{-\int_{0}^{x}\bm{\lambda}(t) dt}\]
from which the proposition follows.
\end{proof}

\subsection{Case of $\bm{\lambda}(t)\simeq \frac{1}{t}$} \label{vn}
When $\bm{\lambda}(t)=\frac{1}{t}$, equation (\ref{eqdiffPhi}) is
the continuous conterpart of equation (\ref{defvn}), so we expect
$u(t)$ to have the same asymptotic behavior as $v_n$. This will be
proved with an additional hypothesis on $\Phi$ in the next section.
Here we show a slightly weaker result without any assumption.

\begin{proposition}
There exists a function $\bm{\lambda}:[0,+\infty]\rightarrow ]0,1]$
such that $\bm{\lambda}(t)\sim\frac{1}{t}$ and for which the
solution $w$ of (\ref{eqdiffPhi}) satisfies
\[\|w(n)-v_n\|\underset{n\rightarrow+\infty}\longrightarrow 0.\]

\end{proposition}

\begin{proof}

Let $U$ be the solution of (\ref{equadiff}) and
$v(t)=\frac{U(t)}{t+1}$, which thus satisfies
\[(t+2)v(t)+(t+1)v'(t)=J((t+1)v(t)).\]
Define $\zeta(t)=t+\ln(1+t)$. By making the change of time
$s=\zeta(t)$ and $w(s)=v(t)$, we get
\[w(s)+w'(s)=\Phi\left(\frac{1}{2+\zeta^{-1}(s)},w(s)\right)\]
and $w$ is thus solution of (\ref{eqdiffPhi}) with
\[\bm{\lambda}(t)=\frac{1}{2+\zeta^{-1}(t)}=\frac{1}{t}+\frac{\ln(t)}{t^2}+o\left(\frac{\ln(t)}{t^2}\right).\]

Moreover, \[\|w(n)-v_n \|\leq
\|v(n)-v_n\|+\left\|v\left(\zeta^{-1}(n)\right)-v(n)\right\|.\] We
already know by Corollary \ref{convvn} that $\|v(n)-v_n\|$ goes to 0
as $n$ tends to $+\infty$. On the other hand, by the mean value
Theorem,

\begin{equation}\label{vzeta}\left\|v\left(\zeta^{-1}(n)\right)-v(n)\right\|\leq
\left(n-\zeta^{-1}(n)\right)\cdot\max_{x\in
[\zeta^{-1}(n),n]}\|v'(x)\|. \end{equation}

By definition of $v$,
$\displaystyle{v'(x)=\frac{(x+1)U'(x)-U(x)}{(x+1)^2}}$ hence
Corollary \ref{decr} implies that

\begin{eqnarray}
\|v'(x)\|&\leq &\frac{(x+1)\|U'(x)\|+\|U(x)\|}{(x+1)^2}\\
&\leq&\frac{(x+1)\|U'(0)\|+\|U(0)\|+x \|U'(0)\|}{(x+1)^2}\\
&\leq&\frac{C}{x+1}
\end{eqnarray}
for $C=2\max(\|U(0)\|,\|U'(0)\|)$.

Replacing in equation (\ref{vzeta}) gives

\[ \left\|v\left(\zeta^{-1}(n)\right)-v(n)\right\|\leq C
\frac{n-\zeta^{-1}(n)}{1+\zeta^{-1}(n)} \] which goes to 0 since
$\zeta(n)\sim n$, and we have thus proved that
\[\|w(n)-v_n\|\underset{n\rightarrow+\infty}\longrightarrow 0.\]
\end{proof}

An interesting corollary of this Proposition, which gives a
sufficient condition for convergence of both $v_n$ and $v_\lambda$
to the same limit, is:

\begin{corollary}\label{corconvboth}
Let $U$ be the solution of (\ref{equadiff}). If $U'(t)$ converges to
$l$ when $t$ goes to $+\infty$, then $v_n$ and $v_\lambda$ converge
to $l$ as well as $n$ goes to $+\infty$ and $\lambda$ goes to 0,
respectively.
\end{corollary}

\begin{proof}
Suppose that $U'(t)$ converges to $l$ when $t$ goes to $+\infty$.
Then $v(t)=\frac{U(t)}{t}$ converges to $l$ as well, and so does
$v_n$ according to Corollary \ref{convvn}.

On the other hand,
\[t v'(t)=U'(t)-\frac{U(t)}{t}\rightarrow l-l=0\]
so $v'(t)=o\left(\frac{1}{t}\right)$. Define $\zeta$, $\bm{\lambda}$
and $w$ as in the proof of the preceding proposition ; then
$w(t)=v(\zeta^{-1}(t))$ converges also to $l$ and by definition

\[w'(\zeta(t))=\frac{t+1}{t+2}\,v'(t)=o\left(\frac{1}{t}\right).\]
Since $\zeta(t)\sim t$ and $\bm{\lambda}(t)\sim\frac{1}{t}$ this
implies that $\frac{\|w'(t)\|}{\bm{\lambda}(t)}=o(1)$. According to
Lemma \ref{simple}, this implies that
$\|w(t)-v_{\bm{\lambda}(t)}\|=o(1)$, and so $v_\lambda$ tends to $l$
as $\lambda$ goes to 0.

\end{proof}

\subsection{Case of a slow parametrization}
%\textrm{}
%\\

From now on the following assumption $(\mathcal{H})$ will be made:
there is a constant $C$ such that
\[\|\Phi(\lambda,x)-\Phi(\mu,x)\|\leq |\lambda-\mu| (C+\|x\|)\quad  \forall x\in X \quad \forall (\lambda,\mu)\in ]0,1]^2.\quad (\mathcal{H})\]
\begin{remark}
$(\mathcal{H})$ is satisfied as soon as $J$ is the Shapley operator
(\ref{Shapley}) of a game with bounded payoff since in that case
\[
\|\Phi(\lambda,x)-\Phi(\mu,x)\|_{\infty}\leq|\lambda-\mu|\left(\|g\|_\infty+\|x\|_\infty\right)
\]

\end{remark}

\begin{remark}\label{deriveevlambda}
Hypothesis $(\mathcal{H})$ implies that for every $\lambda$ and
$\mu$
\[
\frac{\|v_\lambda-v_\mu\|}{|\lambda-\mu|}\leq \frac{C'}{\lambda}
\]
for some constant $C'$: in some sense $(\mathcal{H})$ is thus a
statement about the speed of variation of the family
$\{v_\lambda\}$.
\end{remark}

The principal result of this subsection is Corollary \ref{convder}
which states that under this hypothesis, if the parametrization
$\bm{\lambda}$ converges slowly enough to 0, then the corresponding
solution of (\ref{eqdiffPhi}) has the same asymptotic behavior as
the family $\{v_\lambda\}$. We start by a technical result:

%\begin{proposition}\label{convder}
%If $\lambda$ is $\mathcal{C}^2$, and
%$\frac{\lambda'(t)}{\lambda^2(t)}$ tends to 0 when $t$ goes to
%$+\infty$, then the solution of (\ref{eqdiffPhi}) satisfies
%$\|U(t)-v_{\bm{\lambda}(t)}\|\rightarrow 0$
%\end{proposition}

\begin{proposition}\label{ineqvarlente}
Let $\bm{\lambda}$ be a $\mathcal{C}^1$ function from $[0,+\infty[$
to $]0,1]$ and let $L:\mathbb{R}^+\rightarrow\mathbb{R}$ be defined
by $\displaystyle{L(t)=e^{\int_0^t
\left[\frac{|\bm{\lambda}'(s)|}{\bm{\lambda}(s)}-\bm{\lambda}(s)\right]
ds}}$. Then the corresponding solution $u$ of (\ref{eqdiffPhi})
satisfies:
%\[
%\|u(t)-v_{\bm{\lambda}(t)}\|\leq \frac{e^{\int_0^t
%\left[\frac{|\bm{\lambda}'(s)|}{\bm{\lambda}(s)}-\bm{\lambda}(s)\right]
%ds}}{\bm{\lambda}(t)}\left[ \|u'(0)\|+(C+C')\int_0^t
%|\bm{\lambda}'(s)| e^{\int_0^s
%\left[\bm{\lambda}(r)-\frac{|\bm{\lambda}'(r)|}{\bm{\lambda}(r)}
%\right]dr}ds\right].
%\]
\[
\|u(t)-v_{\bm{\lambda}(t)}\|\leq \frac{L(t)}{\bm{\lambda}(t)}\left[
\|u'(0)\|+(C+C')\int_0^t \frac{|\bm{\lambda}'(s)|}{L(s)} ds\right].
\]
\noindent where $C$ is the constant in condition $(\mathcal{H})$ and
$C'=\sup\limits_{\lambda\in]0,1]} \|v_\lambda\|$.
\end{proposition}

\begin{proof}
For any $h>0$, define $\displaystyle{u_h(t)=\frac{u(t+h)-u(t)}{h}}$
and
$\displaystyle{\bm{\bm{\lambda}}_h(t)=\frac{\bm{\lambda}(t+h)-\bm{\lambda}(t)}{h}}$.
Since $u$ is $\mathcal{C}^1$,
\[
u_h(t)=u'(t)+\frac{1}{h}\int_t^{t+h} u'(s)-u'(t)ds
\]
which implies, by uniform continuity of $u$ on any compact set, that
the restriction of $u_h$ to any closed intervall converges uniformly
to $u'$ as $h$ goes to 0. Similarly, the restriction of
$\bm{\lambda}_h$ to any
closed intervall converges uniformly to $\bm{\lambda}'$ as $h$ goes to 0.\\
Since $u$ satisfies equation (\ref{eqdiffPhi}), for any $h$ and $t$,

\begin{eqnarray}
\|u_h(t)+u_h'(t)\|&=&\frac{1}{h}\|\Phi(\bm{\lambda}(t+h),u(t+h))-\Phi(\bm{\lambda}(t),u(t))\|\\
&\leq&\frac{1}{h}\ \|\Phi(\bm{\lambda}(t+h),u(t+h))-\Phi(\bm{\lambda}(t+h),u(t))\|\\
&&+\nonumber\frac{1}{h}\|\Phi(\bm{\lambda}(t+h),u(t))-\Phi(\bm{\lambda}(t)(t),u(t))\|  \\
&\leq&(1-\bm{\lambda}(t+h))\|u_h(t)\|+|\bm{\lambda}_h(t)|(C+\|u(t)\|).
\end{eqnarray}
\noindent by hypothesis $(\mathcal{H})$. According to Lemma
\ref{simple}, this implies that
\begin{equation}
\|u_h(t)+u_h'(t)\|\leq
(1-\bm{\lambda}(t+h))\|u_h(t)\|+|\bm{\lambda}_h(t)|\left(C+C'+\frac{\|u'(t)\|}{\bm{\lambda}(t)}\right).
\end{equation}
\noindent where $C'$ is a majorant of the family $\|v_\lambda\|$.\\
Fix $t_0>0$, and let $\varepsilon>0$. Since $\bm{\lambda}(t)$ is
bounded from below on $[0,t_0]$ and using the uniform convergence of
$u_h$ to $u'$ on $[0,t_0]$, one obtains that for $h$ small enough,
and for every $t\leq t_0$,
\begin{equation}
\|u_h(t)+u_h'(t)\|\leq
\left(1-\bm{\lambda}(t+h)+\frac{|\bm{\lambda}_h(t)|}{\bm{\lambda}(t)}\right)\|u_h(t)\|+(C+C'+\varepsilon)|\bm{\lambda}_h(t)|.
\end{equation}
Then applying Proposition \ref{decrexpo} to $u_h$ implies that for
any $h$ small enough and $t\leq t_0$,

\[
\|u_h(t)\|\leq\
e^{\int_{0}^{t}\left[\frac{|\bm{\lambda}_h(s)|}{\bm{\lambda}(s)}-\bm{\lambda}(s+h)
\right]ds}\left(\|u_h(0)\|+(C+C'+\varepsilon)\int_0^t
|\bm{\lambda}_h(s)|e^{\int_0^s
\left[\bm{\lambda}(r+h)-\frac{|\bm{\lambda}_h(r)|}{\bm{\lambda}(r)}\right]
dr} ds\right).
\]
Using the uniform convergence of $\bm{\lambda}_h$ and
$\bm{\lambda}(\cdot+h)$ on $[0,t_0]$, letting $h$ go to 0 implies
that for any $t\leq t_0$,
%\[
%\|u'(t)\|\leq
%e^{\int_{0}^{t}\left[\frac{|\bm{\lambda}'(s)|}{\bm{\lambda}(s)}-\bm{\lambda}(s)
%\right]ds}\left(\|u'(0)\|+(C+C'+\varepsilon)\int_0^t
%|\bm{\lambda}'(s)|e^{\int_0^s
%\left[\bm{\lambda}(r)-\frac{|\bm{\lambda}'(r)|}{\bm{\lambda}(r)}
%\right]dr} ds\right).
%\]
\[
\|u'(t)\|\leq L(t)\left(\|u'(0)\|+(C+C'+\varepsilon)\int_0^t
\frac{|\bm{\lambda}'(s)|}{L(s)} ds\right).
\]
Since this is true for every $t_0$ and $\varepsilon$, using Lemma
\ref{simple} again gives
%\[
%\|u(t)-v_{\bm{\lambda}(t)}\|\leq \frac{e^{\int_0^t
%\left[\frac{|\bm{\lambda}'(s)|}{\bm{\lambda}(s)}-\bm{\lambda}(s)
%\right]ds}}{\bm{\lambda}(t)}\left[ \|u'(0)\|+(C+C')\int_0^t
%|\bm{\lambda}'(s)| e^{\int_0^s
%\left[\bm{\lambda}(r)-\frac{|\bm{\lambda}'(r)|}{\bm{\lambda}(r)}
%\right]dr}ds\right].
%\]
\[
\|u(t)-v_{\bm{\lambda}(t)}\|\leq \frac{L(t)}{\bm{\lambda}(t)}\left[
\|u'(0)\|+(C+C')\int_0^t \frac{|\bm{\lambda}'(s)|} {L(s)}ds\right].
\]
\end{proof}

\begin{remark}
If in Proposition \ref{ineqvarlente} we suppose in addition that
$\bm{\lambda}$ is nonincreasing, we get the simpler inequality
\[
\|u(t)-v_{\bm{\lambda}(t)}\|\leq \frac{e^{-\int_0^t \bm{\lambda}(s)
ds}}{\bm{\lambda}^2(t)}\left[ \|u'(0)\|-(C+C')\int_0^t
\bm{\lambda}(s)\bm{\lambda}'(s) e^{\int_0^s
\bm{\lambda}(r)dr}ds\right].
\]
\end{remark}

%An interesting corollary of Proposition \ref{ineqvarlente} holds in
%the case of a slow parametrization:

As a corollary to Proposition \ref{ineqvarlente} we can now prove:

\begin{corollary}\label{convder}
Let $\bm{\lambda}$ be a $\mathcal{C}^1$ function from $[0,+\infty[$
to $]0,1]$, such that $\frac{\bm{\lambda}'(t)}{\bm{\lambda}^2(t)}$
converges to 0 as $t$ goes to $+\infty$, and let $u$ be the
corresponding solution of equation (\ref{eqdiffPhi}). Then
$\|u(t)-v_{\bm{\lambda}(t)}\|$ goes to 0 as $t$ goes to $+\infty$.
\end{corollary}

\begin{proof}
First notice that $\left(\frac{1}{\bm{\lambda}(t)}\right)'=o(1)$, so
$\frac{1}{\bm{\lambda}(t)}=o(t)$ which implies that $\bm{\lambda}(t)\notin\ell^1$.\\
Next we prove that \[\frac{L(t)}{\bm{\lambda}(t)}=o(1).\] Since the
left-hand side is equal to $\displaystyle{\frac{e^{\int_0^t
\left[\frac{|\bm{\lambda}'(s)|}{\bm{\lambda}(s)}-\frac{\bm{\lambda}'(s)}{\bm{\lambda}(s)}-\bm{\lambda}(s)
\right]ds}}{\bm{\lambda}(0)}}$, the result is deduced from the fact
that $\frac{\bm{\lambda}'(s)}{\bm{\lambda}(s)}=o(\bm{\lambda}(s))$
and that $\bm{\lambda}(t)\notin\ell^1$.

 Finally we prove that \[ \int_0^t
\frac{|\bm{\lambda}'(s)|}{L(s)} ds=o\left(
\frac{\bm{\lambda}(t)}{L(t)}\right).
\]
Since the right-hand side diverges to $+\infty$, it is enough to
prove that the derivative satisfies
%\[
%|\bm{\lambda}'(t)| e^{\int_0^t
%\left[\bm{\lambda}(s)-\frac{|\bm{\lambda}'(s)|}{\bm{\lambda}(s)}
%\right]ds}=o\left((\bm{\lambda}'(t)+\bm{\lambda}^2(t)-|\bm{\lambda}'(t)|)e^{\int_0^t
%\left[\bm{\lambda}(s)-\frac{|\bm{\lambda}'(s)|}{\bm{\lambda}(s)}
%\right]ds}\right)
%\]
\[
\frac{|\bm{\lambda}'(t)|}{L(t)}
=o\left(\frac{\bm{\lambda}'(t)+\bm{\lambda}^2(t)-|\bm{\lambda}'(t)|}{L(t)}\right)
\]
\noindent which is true since
$\bm{\lambda}'(t)=o(\bm{\lambda}^2(t))$.
\end{proof}

\begin{remark}
Note the similarity of this proposition with some approximation
results for dynamical systems in the framework of Hilbert spaces,
for example the slow parametrization in \cite{AtCo}:
\begin{itemize}
 \item first there is
a parallel between the strong monotonicity condition in \cite{AtCo}
p. 523 and our assumption that the $\Phi(\lambda,\cdot)$ are
contracting.

\item Second between a condition about the derivative of the
trajectory in the same paper p. 528 and our hypothesis
$(\mathcal{H})$ (see remark \ref{deriveevlambda}).

\item Third the slow-convergence condition is the same (see condition
(ii) in \cite{AtCo} p. 528).

\item Lastly, results of both papers are of the same nature: convergence of a certain family ($\{v_\lambda\}$ in this paper) implies that the solution of any slowly-perturbed
evolution equation tends to this limit as time goes to infinity.
\end{itemize}
 A difference however is the fact that in this paper we also have a
 reciprocal: if for any  slow parametrization $\bm{\lambda}$ the solution $u(t)$ of (\ref{eqdiffPhi})
 converges as $t$ goes to infinity, then the family $v_\lambda$ converges to the same limit as $\lambda$ goes to
 0.
\end{remark}

\begin{remark}
In the proof of Proposition \ref{ineqvarlente} only the three
following hypotheses on the family $\Phi$ were used:
\begin{itemize}
\item[$\mathrm{(i)}$] $\Phi(\cdot,x)$ satisfies condition $\mathcal{H}$ for all $x$.
\item[$\mathrm{(ii)}$] $\Phi(\lambda,\cdot)$ is $1-\lambda$ contracting for every $\lambda\in]0,1]$.
\item[$\mathrm{(iii)}$] The fixed points $v_\lambda$ are uniformly bounded.
\end{itemize}
The two last ones are satisfied as soon as $\Phi(\lambda,x)=\lambda
J\left(\frac{1-\lambda}{\lambda}x\right)$ for a nonexpansive
operator $J$, but this is not a necessary condition for Proposition
\ref{ineqvarlente} to holds.
\end{remark}

\begin{remark}
In fact, the more general result holds: suppose that the family
$\Phi$ satisfies the three hypotheses:
\begin{itemize}
\item[$\mathrm{(i)}$] There exists a constant $C$ and a continuous function $M$ from
$]0,1]$ to $\mathbb{R}^+$ such that for any $(x,\lambda,\mu)$ in
$X\times]0,1]^2$,
\[\|\Phi(\lambda,x)-\Phi(\mu,x)\|\leq
\left|\int_\lambda^\mu M(\gamma)d\gamma\right|(C+\|x\|).\]
\item[$\mathrm{(ii)}$] There exists a continuous function $\beta:]0,1]\rightarrow]0,1]$ such that $\Phi(\lambda,\cdot)$ is $1-\beta(\lambda)$ contracting.
\item[$\mathrm{(iii)}$] The fixed points $v_\lambda$ of $\Phi(\lambda,\cdot)$ are uniformly bounded by $C'$.
\end{itemize}
%Suppose that the $\mathcal{C}^1$
%parametrization$\lambda:[0,+\infty[\rightarrow]0,1]$ satisfies the
%three
%hypotheses:
%\begin{itemize}
%\item $\bm{\lambda}(t)$ goes to 0 in $+\infty$.
%\item $\displaystyle{\liminf_{t\rightarrow +\infty}}\, t\beta(\bm{\lambda}(t)) >1$.
%\item $\frac{\lambda'(t) M(\bm{\lambda}(t))}{\beta^2(\bm{\lambda}(t))}=o(1)$.
%\end{itemize}
Let $\bm{\lambda}$ be a $\mathcal{C}^1$ function from $[0,+\infty[$
to $]0,1]$. Then the corresponding solution $u$ of (\ref{eqdiffPhi})
satisfies \begin{eqnarray*} &&\|u(t)-v_{\bm{\lambda}(t)}\|\leq\\
\,&&\\ &&\frac{e^{\int_0^t
\frac{|\bm{\lambda}'(s)|M(\bm{\lambda}(s))}{\beta(\bm{\lambda}(s))}-\beta(\bm{\lambda}(s))
ds}}{\beta(\bm{\lambda}(t))}\left[\|u'(0)\|+(C+C')\int_0^t
|\bm{\lambda}'(s)|M(\bm{\lambda}(s)) e^{\int_0^s
\left[\beta(\bm{\lambda}(i))-\frac{|\bm{\lambda}'(i)|M(\bm{\lambda}(i))}{\beta(\bm{\lambda}(i))}
\right]di}ds\right]\hspace{-0.8ex}.
\end{eqnarray*}
This implies that $\|u(t)-v_{\bm{\lambda}(t)}\|$ tends to 0 as soon
as $\beta$ is $\mathcal{C}^1$ and the parametrization $\lambda$
satisfies both properties :
\begin{itemize}
\item[$\mathrm{(iv)}$] $\frac{\bm{\lambda}'(t)M(\bm{\lambda}(t))}{\beta^2(\bm{\lambda}(t))}=o(1)$
\item[$\mathrm{(v)}$] $\frac{\bm{\lambda}'(t)\beta'(\bm{\lambda}(t))}{\beta^2(\bm{\lambda}(t))}=o(1)$
\end{itemize}
Notice again the similarity with \cite{AtCo}.

\end{remark}

Another interesting consequence of hypothesis $(\mathcal{H})$ is
Corollary \ref{equiv} which states that if two parametrizations are
close to one other, then this is also the case for the trajectories.
We first prove a technical result using the same approach as in the
proof of Proposition \ref{ineqvarlente}:

\begin{proposition}\label{diffsol}
Let $u$ and $v$ be the two solutions of (\ref{eqdiffPhi}) for some
functions $\bm{\lambda}$ and $\bm{\mu}$ respectively. Then for any
$t\geq 0$,
\[
\|u(t)-v(t)\|\leq e^{-\int_{0}^{t}\bm{\mu}(s)
ds}\left(\|u_0-v_0\|+\int_0^t
\left(C+\|u(s)\|\right)|\bm{\bm{\lambda}}(s)-\bm{\mu}(s)|\cdot
e^{\int_0^s \bm{\mu}(i) di} ds\right)
\]

\end{proposition}

%\begin{proposition} \label{equiv}
%Let $u$ and $v$ the two solutions de (\ref{eqdiffPhi}) for some
%functions $\bm{\lambda}$ and $\bm{\mu}$ respectively. Assume $\bm{\lambda}\sim\bm{\mu}$,
%$\int_0^{+\infty} \bm{\bm{\lambda}}(t)dt=+\infty$ and $u$ bounded, then
%$\|u(t)-v(t)\|\rightarrow 0$.
%\end{proposition}

\begin{proof}
Let $f=u-v$, then
\begin{eqnarray*}\|f(t)+f'(t)\|&=&\|\Phi(\bm{\bm{\lambda}}(t),u(t))-\Phi(\bm{\mu}(t),v(t))\|\\
&\leq&\|\Phi(\bm{\bm{\lambda}}(t),u(t))-\Phi(\bm{\mu}(t),u(t))\| + \|\Phi(\bm{\mu}(t),u(t))-\Phi(\bm{\mu}(t),v(t))\| \\
&\leq&|\bm{\bm{\lambda}}(t)-\bm{\mu}(t)|\cdot(C+\|u(t\|)+(1-\bm{\mu}(t))\cdot\|f(t)\|
\end{eqnarray*}
because of hypothesis $(\mathcal{H})$ and contraction of
$\Phi(\lambda,\cdot)$. Applying Proposition \ref{decrexpo} gives the
result.
\end{proof}

In particular one has:
\begin{corollary}\label{equiv}
Let $u$ and $v$ the two solutions of (\ref{eqdiffPhi}) for some
functions $\bm{\lambda}$ and $\bm{\mu}$ respectively. Assume that
$u$ is bounded and $\bm{\mu}\notin \ell^1$, then
$\|u(t)-v(t)\|\rightarrow 0$ in the two following cases:
\begin{enumerate}
\item[$\mathrm{a)}$] $\bm{\mu}(t)\sim\bm{\bm{\lambda}}(t)$  as $t$ goes to $+\infty$
\item[$\mathrm{b)}$] $|\bm{\lambda}-\bm{\mu}|\in\ell^1.$
\end{enumerate}
\end{corollary}

\begin{proof}
Let $L$ be a bound for $u$. By the preceding proposition we know
that
\[
\|u(t)-v(t)\|\leq e^{-\int_{0}^{t}\bm{\mu}(s)
ds}\left(\|u_0-v_0\|+(C+L)\int_0^t
|\bm{\bm{\lambda}}(s)-\bm{\mu}(s)|\cdot e^{\int_0^s \bm{\mu}(i) di}
ds\right)
\]
so it suffices to show that
\[
\int_0^t |\bm{\bm{\lambda}}(s)-\bm{\mu}(s)|\cdot e^{\int_0^s
\bm{\mu}(i) di} ds=o\left( e^{\int_{0}^{t}\bm{\mu}(s) ds} \right).
\]
\begin{enumerate}
\item[$\mathrm{a)}$] Assume that $\bm{\mu}(t)\sim\bm{\bm{\lambda}}(t)$, that is
$\frac{|\bm{\bm{\lambda}}(t)-\bm{\mu}(t)|}{\bm{\mu}(t)}=o(1)$. This
implies that
\[|\bm{\bm{\lambda}}(t)-\bm{\mu}(t)|\cdot e^{\int_{0}^{t}\bm{\mu}(s) ds}=o\left(\bm{\mu}(t)
e^{\int_{0}^{t}\bm{\mu}(s) ds}\right)\] which gives the result by
integrating.
\item[$\mathrm{b)}$] Assume that $|\bm{\lambda}-\bm{\mu}|\in\ell^1$, let $\varepsilon>0$ and
$T$ such that $\int_T^{+\infty}|\bm{\lambda}(s)-\bm{\mu}(s)|ds
<\varepsilon$. Then for $t>T$,

\begin{eqnarray*}
I&:=&\int_0^t |\bm{\lambda}(s)-\bm{\mu}(s)|\cdot e^{\int_0^s
\bm{\mu}(i) di} ds
\\&=&\int_0^T |\bm{\lambda}(s)-\bm{\mu}(s)|\cdot e^{\int_0^s \bm{\mu}(i) di}
ds+\int_T^t |\bm{\lambda}(s)-\bm{\mu}(s)|\cdot e^{\int_0^s \bm{\mu}(i) di} ds\\
&\leq& e^{\int_0^T \bm{\mu}(s) ds} \int_0^T
|\bm{\lambda}(s)-\bm{\mu}(s)|ds
+e^{\int_0^t \bm{\mu}(s) ds} \int_T^t |\bm{\lambda}(s)-\bm{\mu}(s)|ds\\
&\leq&e^{\int_0^T \bm{\mu}(s) ds} \int_0^T
|\bm{\lambda}(s)-\bm{\mu}(s)|ds+\varepsilon e^{\int_0^t \bm{\mu}(s) ds}\\
&\leq&2\varepsilon e^{\int_0^t \bm{\mu}(s) ds}
\end{eqnarray*}
\noindent for all $t$ large enough since $e^{\int_0^t \bm{\mu}(s)
ds}$ diverges to $+\infty$ as $t$ goes to $+\infty$.
\end{enumerate}
\end{proof}

Some interesting corollaries follows immediately: first because of
Corollary \ref{lambdaconstant}, we get the

\begin{corollary}
If $\bm{\lambda}(t)\rightarrow\lambda>0$, then $u(t)\rightarrow
v_\lambda$
\end{corollary}

%Next, one obtains the refinement of Proposition \ref{convder} :
%
%\begin{proposition}
%Assume $\lambda$ $\mathcal{C}^1$, $\bm{\lambda}'(t)<0$ and that both
%$\bm{\lambda}(t)$ and $\frac{\bm{\lambda}'}{\bm{\lambda}^2}$ tend to 0 as $t$ goes
%to $+\infty$. Then the solution $u$ of (\ref{eqdiffPhi}) satisfies
%$\|u(t)-v_{\lambda(t)}\|\rightarrow 0$
%\end{proposition}
%
%\begin{proof}
%Denote
%$f(t)=\ln(-\frac{\bm{\lambda}'(t)}{\bm{\lambda}^2(t)})=\ln((\frac{1}{\lambda})'(t))$.
%Since $f$ is continuous, one can find a function $g$ $\mathcal{C}^1$
%such that $\|f(t)-g(t)\|\rightarrow 0$ as $t$ goes to $+\infty$.
%This shows that $(\frac{1}{\lambda})'\sim e^g$ in $+\infty$, and
%because $\lambda$ tends to 0 it implies that
%$\frac{1}{\lambda(t)}\sim\int_0^t e^{g(x)}dx$. Define
%$\bm{\mu}(t)=\frac{1}{1+\int_0^t e^{g(x)}dx}$ and notice that $\mu$ is
%$\mathcal{C}^2$, $(\frac{1}{\mu})'=o(1)$ and $\mu\sim\lambda$. Then
%apply Propositions \ref{convder} and \ref{equiv}.
%\end{proof}

Then, combining the results of section \ref{vn} and Corollaries
\ref{convder} and \ref{equiv} we deduce the following Corollary
bringing to light the tight difference between dynamics related to
$\lim v_n$ and $\lim v_\lambda$:

\begin{corollary}\label{coralpha}
For $\alpha\in[0,1[$, let $u^\alpha$ be the solution of
\begin{equation}\label{eqalpha} u(t)+u'(t)=\Phi\left((1+t)^{\alpha-1},u(t)\right) \quad \mathrm{with}\
u(0)=u_0\end{equation} Then $u^0(t)$ converges to some $l\in X$ when
$t$ goes to $+\infty$ iff $v_n$ converges to $l$ as $n$ goes to
$+\infty$ ; and for $\alpha\in]0,1[$ $u^\alpha(t)$ converges to some
$l\in X$ as $t$ goes to $+\infty$ iff $v_\lambda$ converges to $l$
as $\lambda$ goes to 0.
\end{corollary}

\subsection{Back to discrete time}\label{backdiscrete}

We proved in the last section that under hypothesis ($\mathcal{H}$),
the solution of (\ref{eqdiffPhi}) has the same asymptotic behavior
as the family $\left\{v_\lambda\right\}$ as soon as $\bm{\lambda}$
converges slowly enough to 0. One may wonder if it is true as well
in discrete time. For any sequence $(\lambda_n)_{n\in\mathbb{N}}$ in
$]0,1]$, define the discrete counterpart of equation
(\ref{eqdiffPhi}) :
\begin{equation}\label{eqdiscrPhi}
w_n=\Phi(\lambda_n,w_{n-1}) \quad \mathrm{with}\ w(0)=w_0
\end{equation}

Then one obtains the discrete version of Corollary \ref{convder} :
\begin{proposition}\label{discret}
Let $\lambda_n$ be a sequence in $]0,1]$. Assume that both
$\lambda_n$ and $\frac{1}{\lambda_n}-\frac{1}{\lambda_{n+1}}$ tend
to 0 as $n$ goes to $+\infty$. Then the solution
$(w_n)_{n\in\mathbb{N}}$ of (\ref{eqdiscrPhi}) satisfies

\[\|v_{\lambda_n}-w_n\|\rightarrow 0\]
as $n$ goes to $+\infty$.
\end{proposition}

\begin{proof}
The sequence $\gamma_n=\frac{1}{\lambda_n}$ tends to $+\infty$ and
satisfies $\gamma_n-\gamma_{n-1}\rightarrow 0$ as $n$ goes to
$+\infty$. This implies the existence of an interpolation function
$\gamma: \mathbb{R}\rightarrow\mathbb{R}$ which is $\mathcal{C}^2$
and such that for all
 $n$ in $ \mathbb{N}$, $\gamma(n)=\gamma_n$ ,
$\lim_{+\infty}\gamma(t)=+\infty$ and $\lim_{+\infty}\gamma'(t)=0$.
The function $\bm{\lambda}:=\frac{1}{\gamma}$ thus satisfies
$\bm{\lambda}(n)=\lambda_n$ and all the hypotheses of Corollary
\ref{convder}. Let us denote by $u$ the corresponding solution of
equation (\ref{eqdiffPhi}). By Corollary \ref{convder} it is enough
to show that $\|w_n-u(n)\|\rightarrow 0$ as $n$ goes to $+\infty$.

Define $a_n:=\|w_n-u(n)\|$ and let $\varepsilon>0$. Then
\begin{eqnarray}
\notag a_n&=&\|\Phi(\lambda_n,w_{n-1})-\Phi(\lambda_n,u(n))+u'(n)\| \\
\notag &\leq&(1-\lambda_n)\|w_{n-1}-u(n)\|+\|u'(n)\| \\
\notag &\leq&(1-\lambda_n)\|w_{n-1}-u(n-1)\|+\|u(n)-u(n-1)\|+\|u'(n)\| \\
\notag &\leq&(1-\lambda_n)a_{n-1}+2\sup_{t\in[n-1,n]}\|u'(t)\| \\
\notag
&\leq&(1-\lambda_n)a_{n-1}+2\sup_{t\in[n-1,n]}\left\|\frac{u'(t)}{\bm{\lambda}(t)}\right\|
\cdot \sup_{t\in[n-1,n]}\bm{\lambda}(t) \\
\label{eqsn}&\leq&(1-\lambda_n)a_{n-1}+2\varepsilon
\sup_{t\in[n-1,n]}\bm{\lambda}(t)
\end{eqnarray}
for every $n$ large enough because of Corollary \ref{convder}.

Denote $s_n=\displaystyle{\max_{t\in[n-1,n]}}\bm{\lambda}(t)=o(1)$,
and choose $t_n\in[n-1,n]$ such that $s_n=\bm{\lambda}(t_n)$. Let
$T>0$ such that $|\bm{\lambda}'(t)|\leq\bm{\lambda}^2(t)$ for every
$t\geq T$, then by the mean value Theorem, for any $n\geq T+1$,
\begin{eqnarray*}
|s_n-\lambda_n|&=&|\bm{\lambda}(t_n)-\bm{\lambda}(n)|\\
&\leq&|t_n-n|\cdot\sup_{t\in[n-1,n]}|\bm{\lambda}'(t)|\\
&\leq&\sup_{t\in[n-1,n]}\bm{\lambda}^2(t)\\
&=&s^2_n\\
&=&o(s_n)
\end{eqnarray*}
so that $s_n\sim\lambda_n$ as $n$ goes to $+\infty$. Together with
(\ref{eqsn}) this implies that there exists $N$ such that for all
$n\geq N$,
\[a_n\leq(1-\lambda_n)a_{n-1}+3\varepsilon\lambda_n\] and so by induction one prove that for all
$k\in\mathbb{N}$,
\[a_{N+k}-3\varepsilon\leq(
a_{N}-3\varepsilon)\prod_{i=1}^{k}(1-\lambda_{N+i})\] Now
$\frac{1}{\lambda_n}-\frac{1}{\lambda_{n-1}}\rightarrow 0$ implies
that $\frac{1}{n}=o(\lambda_n)$, so the product goes to 0 and we
deduce that $a_{N+k}\leq 4\varepsilon$ for $k$ large enough.
\end{proof}

\begin{corollary}
$v_\lambda$ converges as $\lambda$ goes to 0 if and only if there
exists a sequence $\lambda_n$ satisfying the hypothesis of
Proposition \ref{discret} such that the corresponding sequence $w_n$
defined by (\ref{eqdiscrPhi}) converges.
\end{corollary}
\begin{proof}
Let $\lambda_n$ such that $w_n$ converges. Because of Proposition
\ref{discret}, $v_{\lambda_n}$ converges. Moreover, for all
$\lambda$ and $\mu$, hypothesis $(\mathcal{H})$ implies that,
denoting $C'=\sup\limits_{\lambda\in]0,1]}\|v_\lambda\|$
\begin{eqnarray*}
\|v_\lambda-v_\mu\|&=&\|\Phi(\lambda,v_\lambda)-\Phi(\mu,v_\mu)\|\\
&\leq&\|\Phi(\mu,v_\lambda)-\Phi(\mu,v_\mu)\|+\|\Phi(\lambda,v_\lambda)-\Phi(\mu,v_\lambda)\|\\
&\leq&(1-\mu)\|v_\lambda-v_\mu\|+|\lambda-\mu|(C+C')
\end{eqnarray*}
and thus that
\begin{equation}\label{vlambdamu}\|v_\lambda-v_\mu\|\leq\left|1-\frac{\lambda}{\mu}\right|(C+C').\end{equation}
Since $\lambda_n\rightarrow 0$ and
$\frac{1}{\lambda_n}-\frac{1}{\lambda_{n+1}}\rightarrow 0$,
$|1-\frac{\lambda_{n+1}}{\lambda_{n}}|$ also converges to 0.
Together with inequality (\ref{vlambdamu}) and the fact that
$v_{\lambda_n}$ converges it implies the convergence of $v_\lambda$
as $\lambda$ goes to 0.

Conversely, if $v_\lambda$ converges, then Proposition \ref{discret}
implies that the sequence $w_n$ defined by equation
(\ref{eqdiscrPhi}) converges as soon as $\lambda_n$ and
$\frac{1}{\lambda_n}-\frac{1}{\lambda_{n+1}}$ tend to 0.
\end{proof}

As in the section 3 (Example \ref{uncertain}), there is an
interpretation in terms of games with uncertain duration:
\begin{example}\label{uncertainbis}
Consider the case of a game with Shapley operator $J$. Let
$\{\lambda_n\}$ be a sequence in $]0,1]$ and $w_n$ defined by
equation (\ref{eqdiscrPhi}). Then $w_n$ is the value of the
following game with uncertain duration: with probability $\lambda_n$
the game stops after stage 1, and the payoff is the payoff during
stage 1. With probability $1-\lambda_n$ there is no payoff during
stage 1 but a transition, and game goes to stage 2. Then,
conditionally to the game going to stage 2, with probability
$\lambda_{n-1}$ the game stops after stage 2, and the payoff is the
payoff during stage 2 ; and with probability $1-\lambda_{n-1}$ there
is no payoff during stage 2 but a transition, and game goes to stage
3. If the game goes to stage $n$, with probability $\lambda_1$ the
payoff is the payoff during stage $n$ and with probability
$1-\lambda_1$ the payoff is 0.

\noindent Proposition \ref{discret} then states that if
$\{\lambda_n\}$ is of slow variation, the value of this game with
uncertain duration is close to the value of the
$\lambda_n$-discounted game.
\end{example}

As a final remark to this section, notice the way in which we proved
Proposition \ref{discret}, with a back and forth process to
continuous dynamics ; it should be interesting to search another
proof using only discrete time methods.
\section{Concluding remarks}
\begin{itemize}
\item In this paper we proved that the asymptotic behavior of
$v_n$ and $v_\lambda$ can be derived from the asymptotic behavior of
solutions of some evolutions equations, namely (\ref{equadiff}) and
(\ref{equdiff3}). It should thus be interesting to determine which
additional conditions on the nonexpansive operator $J$ may imply
convergence of the solutions of these equations, and so convergence
of $v_n$ and $v_\lambda$.

\item Notice that Corollary \ref{coralpha} hints that $v_\lambda$ and
$v_n$ should have the same asymptotic behavior for a wide class of
nonexpansive operators, since the study of $\lim v_n$ seems to be a
limit case of the study of $\lim v_\lambda$. Of interest is also
Corollary \ref{corconvboth} which gives a sufficient condition for
existence of both $\lim v_n$ and $\lim v_\lambda$ as well as their
equality.

\item In Examples \ref{uncertain} and \ref{uncertainbis} we saw that
some results that arose naturally during this paper have a nice
interpretation in the framework of games with uncertain duration. In
particular we showed that for specific types of uncertain duration,
the value of those games behave asymptotically either as $v_n$ or
$v_\lambda$ as the expected time played tends to infinity. Following
\cite{Ne,NeSo} it thus should be interesting to study uncertain
duration more generally, hoping that some conditions on the Shapley
Operator will provide convergence of values for more than just
finitely repeated and discounted games.

\section*{Acknowledgments}

This article was written during the course of my PhD thesis. I would
like to thank my advisor Sylvain Sorin as well as Jérôme Bolte, Juan
Peypouquet and an anonymous referee for very helpful comments and
references.
% I would also
%like to thank Roberto Cominetti and Felipe Alvarez.

%\item The main task is now to determine which additional conditions on the nonexpansive operator $J$ may imply convergence of the solutions of (\ref{eqalpha}), and thus of families $v_n$ and $v_\lambda$.
%\item We have not used the fact that Shapley operators defined in (\ref{Shapley}), besides being nonexpansive, also possess the property of preserving the partial order on their domain (see
%\cite{So}).
%\item Notice that Corollary
%\ref{coralpha} hints that $v_\lambda$ and $v_n$ should have the same
%asymptotic behavior for a wide class of nonexpansive operators,
%since the study of $\lim v_n$ seems to be a limit case of the study
%of $\lim v_\lambda$. Of interest is also Corollary \ref{corconvboth}
%which gives a sufficient condition for existence of both $\lim v_n$
%and $\lim v_\lambda$ as well as their equality.
%\item Can we directly prove results of
%Section \ref{backdiscrete} without the passage to continuous time
%dynamics ?

\end{itemize}

\end{document}